\documentstyle[12pt]{amsart}
\newtheorem{theorem}{Theorem}
\newtheorem{lemma}{Lemma}
\newtheorem{proposition}{Proposition}

\newtheorem{definition}{Definition}

\begin{document}
\title{Branch Structure of $J$--holomorphic Curves 
near Periodic Orbits of a Contact Manifold}
\author{Adam Harris and Krzysztof Wysocki}
\date{}
\maketitle

{\bf Abstract:} Let $M$ be a three--dimensional contact 
manifold, and $\tilde{\psi}:D\setminus\{0\}@>>>M\times{\Bbb R}$
a finite--energy pseudoholomorphic map from the punctured disc 
in ${\Bbb C}$ that is asymptotic to a periodic orbit of the 
contact form. This article examines conditions under which
smooth coordinates may be defined in a tubular neighbourhood of
the orbit such that $\tilde{\psi}$ resembles a holomorphic curve,
invoking comparison with the theory of topological linking of
plane complex algebroid curves near a singular point. Examples
of this behaviour which are studied in some detail include 
pseudoholomorphic maps into ${\Bbb E}_{p,q}\times{\Bbb R}$, where
${\Bbb E}_{p,q}$ denotes a rational ellipsoid (contact structure
induced by the standard complex structure on ${\Bbb C}^{2}$), as
well as contact structures arising from non-standard circle--fibrations
of the three--sphere.
     
\vspace{.2in}

{\bf math reviews subject classification:} 32 Q65, \ 53 D10

\newpage

\section{Introduction}

\vspace{.2in}

The local theory of pseudoholomorphic maps from a Riemann surface
into an almost complex, symplectic $4$--manifold $({\cal M}, J, \omega)$
is developed largely around comparisons with the classical theory of
plane algebroid curves. Following the initial investigations of Gromov,
in work of McDuff [12], Micallef and White [14], and Sikorav [15], methods
of construction of local diffeomorphisms between neighbourhoods of 
$p\in{\cal M}^{2k}$ and ${\bold 0}\in{\Bbb C}^{k}$ (specifically, $k=2$)
were found in order to exhibit singular $J$--holomorphic curves 
$\tilde{\psi}:({\Bbb C},0)@>>>({\cal M},p)$ as being locally equivalent
to holomorphic ones. As a result the local topological data associated
with singularities of plane curves can be transferred to the 
pseudoholomorphic context. Recall that the germ of an algebraic curve
$\Gamma$ with singularity at ${\bold 0}\in{\Bbb C}^{2}$ of multiplicity
$n$ is represented by the vanishing locus of a function $f(w_{1},w_{2})
\in{\cal O}_{{\Bbb C}^{2}}({\cal U})$ for some neighbourhood ${\cal U}$ 
of the origin, such that the homogeneous polynomial $H_{n}(w_{1},w_{2})$
of terms of minimal degree $n$ in the Taylor expansion of $f$ at 
${\bold 0}$, has $n\geq 2$. Each linear factor of $H_{n}$ describes a 
complex tangent line to $\Gamma$ at ${\bold 0}$, hence in order that
the germ of $\Gamma$ be irreducible at the origin it is necessary that
$H_{n}$ be the $n$--th power of a single linear term. If ${\Bbb S}^{3}_{
\varepsilon}$ bounds a ball of radius $\varepsilon$ in ${\cal U}$, then 
the geometric locus of this term (i.e., ignoring multiplicity) intersects
${\Bbb S}^{3}_{\varepsilon}$ in a great circle $C$ corresponding to the
axis of a solid torus in which $\Gamma\cap{\Bbb S}^{3}_{\varepsilon}$ also
describes the trajectory of an iterated torus knot $K_{\Gamma}$. The 
topological ``linking'' of $K_{\Gamma}$ with $C$ provides a set of 
numerical invariants in addition to the multiplicity, on which the 
classification of singular curve--germs is based (cf., e.g., [2]).

The existence of a local analytic {\em parametrization} of $\Gamma$ is a
problem first studied systematically by Newton, and later developed 
rigorously by Weierstrass and Puiseux among others. Thus $\Gamma$ may be
represented in standard form as the image of a holomorphic map
$\Phi:(D,0)@>>>({\cal U},{\bold 0}) \ , \ D\subseteq{\Bbb C}$, such that
$\Phi(z) = (z^{n},\Sigma_{i\geq n+1}a_{i}z^{i})$ (cf. [2]). This type 
of series expansion takes advantage of the natural splitting of ${\Bbb
C}^{2}$ into a product of coordinate planes. In order to see conveniently 
the relationship between certain essential exponents of the parametrization
and the linking invariants of $K_{\Gamma}$ most treatments therefore replace
the ball bounded by ${\Bbb S}^{3}_{\varepsilon}$ with a bidisc $\Delta\times
\Delta$ (this is after all homeomorphic to the ball, hence topological data
are preserved). On the other hand, the totality of all complex 
$1$--dimensional subspaces of ${\Bbb C}^{2}$ induces a fibration of ${\Bbb
S}^{3}$ by great circles, discovered by Hopf. Given a vector ${\bold v}_{p}
\in{\Bbb C}^{2}$ corresponding to $p\in{\Bbb S}^{3}$, the standard complex
structure $J_{0}$ defines a vector field $X(p):= J_{0}\cdot{\bold v}_{p}$
which is tangent to the great circle through each $p\in{\Bbb S}^{3}$. Via
the Euclidean inner product on ${\Bbb R}^{4}\approx{\Bbb C}^{2}$, the
orthogonal complement of $X$ in $T{\Bbb S}^{3}$ defines a non--integrable
plane--field, or ``contact structure'', on ${\Bbb S}^{3}$ corresponding to 
the kernel of the $1$--form $\lambda_{0}$ which is the metric dual of $X$.
The splitting of $T{\Bbb C}^{2}\mid_{{\Bbb S}^{3}}$ induced by this 
structure provides an alternative natural frame for the relationship between
$K_{\Gamma}$ and holomorphic parametrization of plane curves, though we will
return to this matter only in the final section as the ``classical model''
of a more general study.

Contact structures are defined on manifolds $M$ of odd dimension, though
for the purposes of this article we will always assume the dimension to be
three. If $\lambda$ is a $1$--form representing such a structure, the 
extension of $\omega:= d\lambda$ to $M\times{\Bbb R}$ determines a 
``symplectisation'' of $(M,\lambda)$. Conversely, compact symplectic manifolds
$({\cal M},\omega)$ with ends of symplectisation type motivate the role of
contact geometry within modern symplectic topology. Through the
work of Hofer, Wysocki and Zehnder [8], [9] and Eliashberg [3],
$J$--holomorphic curves have been adapted via symplectisation to contact
geometry and topology, specifically as a tool for analysing the Weinstein
conjecture [6] and conversely for explicit construction of moduli of 
pseudoholomorphic maps into manifolds ${\cal M}$ with symplectisation--type
ends [3], [9]. Recall that the ``Reeb vector field'' $X$ associated with 
a given structure $\lambda$ is uniquely determined by the conditions $\lambda(X)\equiv
1 \ , \ {\cal L}_{X}\lambda = 0$ (Lie symmetry). The Weinstein conjecture
asserts that the Reeb flow of a compact contact $3$--manifold always admits
at least one closed characteristic (or ``periodic orbit'') $\gamma:{\Bbb S}
^{1}@>>>M$ such that $\dot{\gamma}(t) = X(\gamma(t))$. If $\xi\subset TM$
denotes the plane--field corresponding to $ker(\lambda)$, then a contact
structure on $M$ is said to be ``pseudohermitian'' when it is equipped
with a partial almost complex structure $j\in C^{\infty}
(M,\xi\otimes\xi^{*})$ such that $d\lambda(*,j\cdot *)\mid_{\xi}$ is a 
positive definite and symmetric quadratic form. The splitting of $TM$ via
$X$ and $\xi$ also provides a natural extension of $j$ to $T(M\times{\Bbb R})$
in relation to which the Cauchy--Riemann equation is defined for 
pseudoholomorphic maps into $M\times{\Bbb R}$ (cf., e.g., [6]). 

Let $\tilde{\psi}:\Sigma@>>>M\times{\Bbb R}$ be one such map. 
For any smooth function $h:{\Bbb R}\rightarrow[0,1]$ we may extend the contact
form to $M\times{\Bbb R}$ by defining
\[\lambda_{h}(p,a) = h(a)\cdot\lambda(p)\]
to act on each $T_{(p,a)}(M\times{\Bbb R})$. In particular, if ${\cal F}_{+}^{1}$
denotes the space of all smooth $h$ such that the derivative $h'\geq 0$, then 
the ``energy'' is given as 
\[E(\tilde{\psi}) = \sup_{h\in{\cal F}_{+}^{1}}\int_{\Sigma}\tilde{\psi}^{*}d\lambda_{h} \ .\]
If $\psi$ denotes the projection of $\tilde{\psi}$ onto $M$ and $a$ its
projection onto ${\Bbb R}$, recall that
\[ \psi^{*}d\lambda = (|\psi_{\eta}|^{2} + |\psi_{\zeta}|^{2})d\eta\wedge
d\zeta \ ,\]
where $z = \eta + {\bold i}\zeta$ is a local complex coordinate on $\Sigma$.
$E(\tilde{\psi})$ consequently vanishes if and only if $\tilde{\psi}$ is
constant -- a condition that follows automatically from Stokes' theorem
when $\Sigma$ is compact, hence pseudoholomorphic maps into symplectisations
are naturally defined on {\em punctured} Riemann surfaces. Specifically,
let $\tilde{\psi}$ be a finite--energy map defined on a punctured disc
$D\setminus\{0\}\subset{\Bbb C}$ and let $r = -\ln(|z|)$ and $\varphi = 
\arg(z)\in{\Bbb S}^{1}\approx{\Bbb R}/2\pi{\Bbb Z}$. A theorem of Hofer,
Wysocki and Zehnder [8] states that if $\tilde{\psi}$ has unbounded image
in $M\times{\Bbb R}$ then there exists a real number $T\neq 0$ and a 
Reeb--periodic orbit $\gamma:{\Bbb S}^{1}@>>>M$ such that 
\[\lim_{r@>>>\infty}\psi(r,\varphi) = \gamma(T\varphi)\hspace{.1in};
\hspace{.1in} \lim_{r@>>>\infty}\frac{a(r,\varphi)}{r} = T\hspace{.1in}
\hbox{in}\hspace{.1in} C^{\infty}({\Bbb S}^{1}) \ .\]
$T$ is moreover an integer multiple of the minimal period $\tau$ of $\gamma$
and corresponds to the ``charge'' $\lim_{r@>>>\infty}\frac{1}{2\pi}\int_{{
\Bbb S}^{1}}\psi^{*}\lambda$ of $\tilde{\psi}$ at $0\in D$. Closed 
characteristics of the Reeb flow are thus realised asymptotically by 
cylinders mapped pseudoholomorphically into the corresponding 
symplectisation of $M$ (we pass over the very substantial theory devoted
to {\em existence} of such mappings of cylinders in general; cf., 
however, [6],[8]). 

The asymptotic relations established in this theorem and refined elsewhere
[10] are fundamental to the present article in which we examine the 
topological behaviour of $\tilde{\psi}(D\setminus\{0\})$ within a tubular
neighbourhood of ${\cal P}:=\gamma({\Bbb S}^{1})\subset M$. Before summarising
our method and results, however, it should be mentioned that they require 
added technical hypotheses to be imposed locally upon the Reeb flow itself. 
Let $\Delta\times({\Bbb R}/ \tau{\Bbb Z}) \ , \ (\Delta\subseteq{\Bbb R}^{2})$
be the tubular domain on which a general system of ordinary differential 
equations $\dot{{\bold x}} = f({\bold x})$, for $f$ smooth and 
$\vartheta$--periodic (minimal period $\tau$) in the coordinates 
${\bold x} = (x,y,\vartheta)\in
\Delta\times{\Bbb S}^{1}$, is defined. Suppose that $\{x = y = 0\}$ corresponds
to a periodic solution $\gamma(\vartheta)$, i.e., if $\phi(x,y,\vartheta)$
denotes the associated $\vartheta$--parametric family of diffeomorphisms, 
depending on initial conditions $(x,y)\in\Delta$ suffciently small, then 
$\gamma(\vartheta) = \phi(0,0,\vartheta) \ , \ \gamma(\vartheta+\tau) = 
\gamma(\vartheta)$. The linear variational equation
for $\dot{{\bold x}} = f({\bold x})$ at $\gamma$ then has the form
\[\frac{\partial}{\partial\vartheta}\phi_{*}(0,0,\vartheta) = f_{*}(0,0,
\vartheta)\cdot\phi_{*}(0,0,\vartheta) \hspace{.1in}(\dag),\]
noting that $f_{*}(0,0,0) = f_{*}(0,0,\tau)$ and $\phi_{*}(0,0,\tau) = 
A\cdot\phi_{*}(0,0,0) = A$, where $A$ denotes the ``holonomy matrix''
associated with (\dag). From the theory of ordinary differential equations
(cf. e.g., [5]) it is well--known that the eigenvalues, or ``characteristic
multipliers'' of $A$ are important in determining the stability of the 
flow $\phi$ near an orbit ${\cal P}$. In particular, $1$ is automatically
a characteristic multiplier in the $\vartheta$--direction, while eigenvalues 
of complex modulus greater than or less than one imply that the flow will be
either unstable or asymptotically stable along closed surfaces containing
${\cal P}$ (cf. [5]). Specifically when $f$ corresponds to the Reeb vector
field of a contact structure, the flow $\phi$ is {\em area--preserving}, 
hence $det(A) = 1$ and the remaining pair of eigenvalues are mutually 
reciprocal (mutually conjugate if they are unimodular). Hofer, Wysocki
and Zehnder actually use the asumption that $A$ is not the identity (cf.
``non--degeneracy'' of ${\cal P}$) as a technical hypothesis in obtaining
their asymptotic formulae (cf. [8]). From the viewpoint of the present 
article this condition has the disadvantage that it does not hold for the
Reeb flow of the standard contact structure on the $3$--sphere, corresponding
to the Hopf fibration as mentioned above. In fact, it does not hold for a
large sub--class of the contact structures which are examined below, although 
an extension of their methods does allow these authors to include some
contact manifolds that are foliated by periodic orbits (cf. [10]). 
The asymptotics have been rederived in a recent thesis of F. Bourgeois [1], 
however, by means of the following alternative hypothesis.

\begin{definition} (cf. [1]) A contact form $\lambda$ on $M$ is said to
be of ``Morse--Bott type'' if, for every $T>0 \ , \  N_{T} := \{p\in M \ 
| \ \phi(p,T) = p\}$ is a smooth, closed, orientable submanifold of $M$
such that $d\lambda\mid_{N_{T}}$ has locally constant rank, and $T_{p}
N_{T} = ker(\phi_{*}(p,T) - Id)$.

\end{definition}    

For present purposes it will henceforth be assumed that all contact structures
under consideration are of Morse--Bott type. Our aim is to understand the 
extent to which $J$--holomorphic cylinders in symplectisations actually 
{\em link} with the periodic orbits which they approach asymptotically. In
this sense the case of algebroid knots $K_{\Gamma}$ linking with great
circles in the $3$--sphere suggests a classical model, since ${\Bbb S}^{3}
\times{\Bbb R}$ is clearly diffeomorphic to ${\Bbb C}^{2}\setminus\{0\}$
(we will return to this matter below). By analogy with the work of Sikorav
and others we will show that a suitably constructed diffeomorphism applied
to a tubular neighbourhood of ${\cal P}$ allows $\tilde{\psi}$ to be 
represented locally by a holomorphic parametrization, and the sense in which
this is possible will be the subject of sections two and three. The problem
of ``asymptotic similarity'' addressed here has a global aspect, however, not
found in the study of singularities of $J$--holomorphic curves. Because of
this we have imposed two further conditions on the pseudohermitian 
structure :

$(1)$ in a sense to be made precise in section three, ${\cal P}$ will be said
to be ``locally recurrent'' if a sufficiently small disc transverse to ${\cal
P}$ can be considered a surface of section of the Reeb flow (the 
diffeomorphism of the disc consequently induced by the return map will be 
denoted by $\alpha$). In particular this condition forces the characteristic
multipliers of $A$ to be unimodular, but is a stronger assumption than orbital
stability near ${\cal P}$; 

$(2)$ it will be assumed that the Reeb flow is not only a Lie symmetry
of the contact form, but also the almost complex structure, i.e., 
${\cal L}_{X}j = 0$ (hence in particular $\alpha$ is diffeomorphically 
equivalent to a rotation). 

It remains now to summarise our results.                 

\vspace{.1in}

{\bf Theorem 1}(cf. section 3) {\em Let $(\psi,a):D\setminus\{0\}@>>>M
\times{\Bbb R}$ be a $J$--holomorphic curve of finite energy and 
charge $n$ at $z=0$,
asymptotic to a locally recurrent periodic orbit ${\cal P}$, near which 
${\cal L}_{X_{\lambda}}j = 0$. Consider any tubular neighbourhood
of ${\cal P}$ in $M$, diffeomorphic to $\Delta\times{\Bbb S}^{1}$ such
that $\{0\}\times{\Bbb S}^{1}\approx{\cal P}$ and such that $\psi^{-1}(
\Delta\times\{\vartheta_{0}\})$ for some fixed $\vartheta_{0}\in{\Bbb S}^{1}$
divides $D\setminus\{0\}$ into $n$ ``quasi--sectors'' $Q_{k}$. 
Then there exists a diffeomorphic
change of coordinates in $\Delta\times[0,2\pi)$ such that on each 
$Q_{k}\subset D\setminus\{0\}$ the map $(\psi,a)$ 
can be expressed in the form 
\[(F_{k}(z), H_{k}(z) - \frac{1}{2\pi{\bold i}}\hat{G}(z,\bar{z}))\ , \ 
\hspace{.1in} 0\leq k\leq n-1 ,\]
where $F_{k} \ , \ H_{k}$ are holomorphic on $Q_{k}$ and continuous
on $\overline{Q}_{k}$, such that
\[F_{k}\mid_{\overline{Q}_{k}\cap\overline{Q}_{k+1}} = \alpha\circ 
F_{k+1}\mid_{\overline{Q}_{k}\cap\overline{Q}_{k+1}} \ , \]
while each $H_{k}$ corresponds to an
analytic branch of $\frac{1}{2\pi{\bold i}}\log(\rho)$ \ , for $\rho$
holomorphic such that $ord_{0}(\rho) = n$.
Moreover, the function $\hat{G}$ belongs to $C^{2}(D)$ and is bounded 
by $K|z|$. Finally, if $\alpha = 1$, then each $F_{k}$ is 
the restriction of a single function $F$ holomorphic on $D$, \ $F(0) = 0$}.

\vspace{.1in}

A representation corresponding to the classical local parametrization
of plane algebroid curves is moreover easily derived when $\alpha = 1$. 
In this case the relationship between  charge and algebraic
multiplicity also becomes explicit. 

\vspace{.1in}

In section four we construct a class of pseudohermitian structures
on the $3$--sphere for which $\alpha = 1$, containing the standard structure
as a distinct special case, and using results on characterisation of 
generalised circle fibrations due to Gluck and Warner [4]. Recent ideas of
McKay [13] relating elliptic line congruences to ``osculating'' almost 
complex structures associated with a four--dimensional real vector space
$V$ are also incorporated in the construction.

\vspace{.1in}

{\bf Theorem 2}(cf. section 4) {\em Let $J$ be the osculating complex 
structure of an elliptic line congruence $\Sigma\subset{\Bbb G}r_{2}(V)$, 
such that the associated 
skew--symmetric 2--form $\omega_{\Sigma}$ on $V\setminus\{0\}$
is closed. Then the 1--form $\lambda$, such that $\lambda_{{\bold v}}
:=i_{{\bold v}}\omega_{\Sigma}$, defines a (fillable, hence tight)
contact structure on ${\Bbb S}^{3}$ for which ${\cal L}_{X_{\lambda}}
J\mid_{ker(\lambda)} = 0$ and $\alpha = 1$. Moreover, two such structures 
are equivalent via a diffeomorphism $\delta$ of ${\Bbb S}^{3}$ if and only 
if $\delta\in{\Bbb O}(4)$}.

\vspace{.1in}

(The notions of ``tight'' and ``over--twisted'' contact structures are
important in the study of global dynamics and topology (cf. [3],[6])
but are not referred to explicitly in this article.)

The final section is devoted to two examples, the first of which is the 
standard structure on ${\Bbb S}^{3}$. Holomorphic maps $\Phi:(D,0)@>>>
({\Bbb C}^{2},{\bold 0})$ corresponding to irreducible germs of a curve
$\Gamma$ are translated explicitly into pseudoholomorphic maps 
$\Psi:D@>>>{\Bbb S}^{3}\times{\Bbb R}$ such that the tangent to $\Gamma$
at the origin becomes the locally recurrent periodic orbit ${\cal P}$ 
corresponding to a great circle. We note that $\alpha = 1$ in this case as a
direct consequence of the Hopf fibration restricted to any tubular 
neighbourhood of ${\cal P}$. It should also be mentioned that the 
identification of algebraic curves with finite--energy pseudoholomorphic
maps into $({\Bbb S}^{3}\times{\Bbb R}, \lambda_{0}, J_{0})$ has been 
addressed from a different 
perspective in an article of Hofer and Kriener [7]. The second example 
discussed in this section includes the well--known contact structures 
$\lambda_{0}$ more generally induced by $J_{0}$ and the Euclidean inner 
product on rational ellipsoids
\[{\Bbb E}_{p,q} = \{(w_{1},w_{2})\in{\Bbb C}^{2} \ | \ p|w_{1}|^{2} +
q|w_{2}|^{2} = 1\} \ ,\]
where $\frac{p}{q}\in{\Bbb Q}$. In this case the hypersurface is foliated
by recurrent periodic orbits such that either $\alpha$ is the identity
or a non--trivial rational rotation through either $2\pi\frac{p}{q}$ or
$2\pi\frac{q}{p}$. Given a tubular
neighbourhood $\Delta\times{\Bbb S}^{1}$ of one of the two periodic orbits
${\cal P}$ for which $\alpha$ is non--trivial, we construct a solid torus 
$\Delta'\times{\Bbb S}^{1}$ with contact form $\lambda$, and a smooth 
covering map $\beta:\Delta\times{\Bbb S}^{1}@>>>\Delta'\times{\Bbb S}^{1}$
such that $\lambda_{0} = \beta^{*}\lambda$ and $\beta_{*}\circ j_{0} = 
j_{0}'\circ\beta_{*}$ (where $j_{0} = J_{0}\mid_{ker(\lambda_{0})}$ and
$j_{0}' = J_{0}\mid_{ker(\lambda)}$). Pseudoholomorphic maps asymptotic
to this orbit in $\Delta
\times{\Bbb S}^{1}$ then project pseudoholomorphically via $\beta$ into
$\Delta'\times{\Bbb S}^{1}$, where the central axis is a recurrent orbit
of the projected Reeb vector field such that $\alpha' = 1$. Working backwards 
from a holomorphic parametrization in $\Delta'\times{\Bbb S}^{1}\times
{\Bbb R}$, we then explicitly construct pseudoholomorphic maps 
$(\psi,a):D\setminus\{0\}@>>>\Delta\times{\Bbb S}^{1}\times{\Bbb R}$ 
according to the prescription of Theorem 1 (cf. section 5).

The authors would like to thank Dr S. Gadde and Dr Y. Tonegawa for 
helpful and stimulating discussions on occasions during the period
of research for this article. The author (Harris) also gratefully
acknowledges the support of an Australian Research Council grant
during the initial phase of work undertaken at the University of 
Melbourne.
 
\vspace{.2in}

\section{Local Model of parametrization}

\vspace{.2in}

In this section we will first review the phenomenon of holomorphic 
similarity in the neighbourhood of a smoothly embedded point, using 
this as a model for the asymptotic case to follow.
Let $M$ be a closed, oriented three--manifold with contact form
$\lambda$ and associated Reeb vector field $X_{\lambda}$. If 
$\xi\subset TM$ denotes the sub-bundle corresponding to the
kernel of $\lambda$, let $J$ be an almost complex structure on
$\xi$ which is compatible with $\lambda$. Consider a Riemann surface
$\Sigma$ and a pseudoholomorphic map
\[(\psi,a):\Sigma@>>>M\times{\Bbb R}.\]
Let $p_{0}\in im(\psi)$ be a smoothly embedded point, and let $D$ be a
neighbourhood of $\psi^{-1}(p_{0})$ with complex coordinate $z = \eta +
{\bold i}\zeta$ such that $z(\psi^{-1}(p_{0})) = 0$. 
For all $p\in M$, and all $v\in T_{p}M$, consider the projection map
$\pi_{p}:T_{p}M@>>>\xi_{p}$ 
defined by $\pi_{p}(v) = v - \lambda_{p}(v)\cdot X_{\lambda}$.
$\psi$ then satisfies the equations 
\[\pi(\psi_{\eta}) + J\pi(\psi_{\zeta}) = 0 \ \hspace{.1in}(*) \ ,\]
\[\lambda(\psi_{\zeta}) = -a_{\eta} \ \hspace{.1in} \ 
\lambda(\psi_{\eta}) = a_{\zeta} \ \hspace{.1in}(\dag) \ .\]
We recall as a consequence that
\[\psi^{*}d\lambda = (|\pi(\psi_{\eta})|^{2} + 
|\pi(\psi_{\zeta})|^{2})d\eta\wedge d\zeta \ ,\]
hence $X_{\lambda}(p_{0})$ does not belong to the image of 
$\psi_{*}$ at the origin of $D$ if and only if 
$|\pi(\psi_{\eta})|^{2} + |\pi(\psi_{\zeta})|^{2}$
does not vanish there  -- a fact which is in turn granted by the 
assumption that $p_{0}$ is smoothly embedded.  Now choose 
a system of Darboux coordinates $(x_{1},x_{2},x_{3})$ neighbouring
$p_{0}$ in $M$. Within this neighbourhood the image of $\psi$ may
then be realized as the graph of a smooth function $x_{3} = 
f(x_{1},x_{2})$. Noting that $\lambda = dx_{3}+x_{1}dx_{2}$ and 
$X_{\lambda} = \frac{\partial}{\partial x_{3}}$
in these coordinates, we remark that 
\[\pi_{p}(v_{1},v_{2},v_{3}) = (v_{1},v_{2},-x_{1}v_{2}) \ ,\]
hence the standard projection
$(v_{1},v_{2},v_{3})\mapsto (v_{1},v_{2})$ defines a linear isomorphism
$\mu$ between $\xi_{p}$ and ${\Bbb R}^{2}$ for all $p$ near the origin at
$p_{0}$. It follows that $\mu\circ J\circ\mu^{-1}$ corresponds to a
$2\times 2$ matrix $j$, the restriction of which to the graph of 
$im(\psi)$ depends only on the coordinates ${\bold x}:= (x_{1},x_{2})$. 
Under projection by $\mu$, equation (*) then becomes 
\[{\bold x}_{\eta}(\eta,\zeta) + j({\bold x})\cdot{\bold x}_{\zeta}
(\eta,\zeta) = {\bold 0} \ .\]
Letting ${\bold e}_{1}$ denote the vector $(1,0)$, consider the system
of ordinary differential equations 
\[\frac{d{\bold x}}{dt} = j({\bold x})\cdot{\bold e}_{1} \ ,\]
and define a local diffeomorphism $\varphi$ near $p_{0}$, via existence 
and uniqueness of solutions, such that
\[(\dot{x}_{1},\dot{x}_{2}) = j(x_{1},x_{2})\cdot{\bold e}_{1} \ 
\Rightarrow\varphi(x,t) := (x_{1}(t),x_{2}(t)) \ , \ \hspace{.2in}
\varphi(x,0) = (x,0) \ .\]
Hence 
\[\varphi_{*}^{-1}\circ j\circ\varphi_{*} = j_{0} = 
\left(\begin{array}{cc}
0 & -1 \\
1 & 0
\end{array}\right) \ , \]
and $\tilde{\bold x}_{\eta}:=\varphi_{*}^{-1}({\bold x}_{\eta}) \ , \ $ 
$\tilde{\bold x}_{\zeta}:=\varphi_{*}^{-1}({\bold x}_{\zeta})$ implies
\[\tilde{\bold x}_{\eta}(\eta,\zeta) + j_{0}\cdot
\tilde{\bold x}_{\zeta}(\eta,\zeta) = {\bold 0} \ .\]
Or, more simply, writing $w = x+{\bold i}t$, we have $w = F(z)$
on a possibly smaller disc $D'\subseteq D\subset{\Bbb C}$, containing
the origin, such that $F(0) = 0 \ , \  \frac{\partial F}
{\partial\bar{z}} = 0$.

\vspace{.1in}

Within the coordinate neighbourhood defined by 
$(x, t, x_{3})$, we now revisit the equations $(\dag)$, with respect
to 
\[\lambda' = (\varphi\times 1)^{*}\lambda = dx_{3}+f_{1}(x,t)dx + 
f_{2}(x,t)dt \ ,\] 
in the form
\[\lambda'(\Psi_{\eta}) = (x_{3})_{\eta} + f_{1}(x,t)x_{\eta} 
+ f_{2}(x,t)t_{\eta} = a_{\zeta}\]
\[\lambda'(\Psi_{\zeta}) = (x_{3})_{\zeta} + f_{1}(x,t)x_{\zeta} 
+ f_{2}(x,t)t_{\zeta} = -a_{\eta} \ ,\]
where $\Psi:=(\varphi\times 1)^{-1}\psi$.
Hence $u = x_{3}+{\bold i}a$ implies
\[\frac{\partial u}{\partial\bar{z}} = -\{f_{1}(x,t)\frac{\partial x}
{\partial\bar{z}} + f_{2}(x,t)\frac{\partial t}{\partial\bar{z}}\}\]
\[ = -\{f_{1}(F(z),\overline{F}(z))\frac{\partial}{\partial\bar{z}}
\Re(F) + f_{2}(F(z),\overline{F}(z))\frac{\partial}{\partial\bar{z}}
\Im(F)\}\]
\[ = -\frac{1}{2}[(f_{1}+{\bold i}f_{2})\circ F(z)]\cdot\overline{
F'(z)} := -F^{*}\omega \ ,\]
where $\omega :=
\frac{1}{2}(f_{1}(w,\bar{w})+{\bold i}f_{2}(w,\bar{w}))d\bar{w}$.
Let $\Delta$ denote a small disc centred at the origin of the 
$(x,t)$--plane, and apply the Cauchy--Green formula to obtain a 
smooth function $g(w,\bar{w}) \ , \ g(0) = 0$ satisfying the equation 
$\bar{\partial}g = \omega$ on $\Delta$. Then we have
\[F^{*}\omega = \bar{\partial}(F^{*}g) = \bar{\partial}(g\circ F) \ .\]
Hence for the parametrization of $\Psi$ we now have the system of equations
\[w = F(z) \ ; \ \hspace{.1in} u = H(z) - g\circ F(z) \ ,\]
with respect to holomorphic functions $F$ and $H$, which simplifies 
under the coordinate transformation 
\[w' = w \ ; \ \hspace{.1in} u' = u + g(w,\bar{w}) - {\bold i}a(0) \]
to $u' = H(z) - {\bold i}a(0)$. Finally, 
$ord_{0}(H -{\bold i}a(0)) = n$ implies
there exists a holomorphic function $h(z) \ , \ h(0) = 0 \ ,
 \ h'(0)\neq 0$, on a possibly smaller disc $D'\subseteq D$,  such that
$H -{\bold i}a(0) = h^{n}$ \ , \ hence $z':=h(z)$ implies 
\[w' = F\circ h^{-1}(z') \ ; \ \hspace{.1in} u' = (z')^{n} \ ,\]
which corresponds to the classical local parametrization of algebroid
curves in ${\Bbb C}^{2}$. In the following section we will examine
a class of periodic orbits near which the analytic 
representation of the local model above can be achieved in a similar 
manner.

\section{Asymptotic approximation near a periodic orbit}

\vspace{.1in}

Consider a periodic orbit of the Reeb flow, denoted ${\cal P}$,
and a tubular neighbourhood $T_{\cal P}\subset M$.
If $\Delta$ represents a disc centred at the origin 
in ${\Bbb R}^{2}$, let $\tilde{\Delta}\subset M$ be an embedded image
such that the origin is mapped to the unique element $p_{0}$
of ${\cal P}\cap\tilde{\Delta}$, with $\tilde{\Delta}$ itself 
corresponding to a 
transverse slice of $T_{\cal P}$. The Reeb flow will be assumed moreover
to be Lyapunov--stable near ${\cal P}$ in the sense that for all 
$p\in\tilde{\Delta'}$,
where $\Delta'\subseteq\Delta$ is a sufficiently small disc centred at
the origin, there exists a unique solution 
$\gamma_{p}:[0,\infty)@>>>M$ to the equation
\[\frac{d\gamma_{p}}{dt} = X_{\lambda}(\gamma_{p}(t)) \ , \ \hspace{.2in}
\gamma_{p}(0) = p \ ,\]
which depends smoothly on both $t$ and $p$, and remains inside $T_{\cal P}$
for all $t\geq 0$. Given $p\in\tilde{\Delta'}$, we 
will define $(i) 
\ \tau(p)$ to be the smallest $t>0$ such that $\gamma_{p}(t)\in
\tilde{\Delta} \ , \ (ii) \ \Gamma_{p} := \gamma_{p}((0,\tau(p)])$
 \ and for each connected open neighbourhood of the origin $\Omega
\subseteq\tilde{\Delta}$ , 
\[(iii) \ \Gamma(\Omega) := \cup_{p\in\Omega\times\{\theta_{0}\}}
\Gamma_{p} \ .\]
We may now consider a recursively defined system of neighbourhoods
$\{\Omega_{k}\}$, such that $\Omega_{0}:=\tilde{\Delta'}$, while $\Omega_{k}$
denotes the origin--component of $\Gamma(\Omega_{k-1})\cap\Omega_{k-1}$.
The set $\Omega_{\infty}:=\cap_{k=0}^{\infty}\Omega_{k}$ measures an
important dynamic aspect of the Reeb flow. 

\begin{proposition} If $\Omega_{\infty}$ is open, then it is 
conformally equivalent to a disc.
\end{proposition} 

\begin{pf} Conformal equivalence to a disc will follow immediately
if $\Omega_{\infty}$ can be shown to be simply connected, which we now
prove by induction : $\Omega_{0} = \tilde{\Delta'}$. Suppose $\Omega_{k}$
is simply connected, and let $C\subset\Omega_{k+1}$ be a simple closed
loop. By the existence and uniqueness of ordinary differential equations
$C$ can be traced back under the Reeb flow to a simple closed loop $C'$ 
which bounds a contractible subdomain, say ${\cal U}\subset\Omega_{k}$.
Once again, existence and uniqueness ensures that for all $p\in{\cal U}$
we have $\Gamma_{p}\cap\Gamma(C') = \emptyset$, while the Jordan Curve
Theorem for plane domains implies that the domain bounded by a given closed
simple loop is unique. It follows that
\[C = \Gamma(\partial{\cal U})\cap\Omega_{k} = \partial\Gamma({\cal U})
\cap\Omega_{k} = \partial(\Gamma({\cal U})\cap\Omega_{k}) \ .\]
Moreover, $\Gamma({\cal U})\cap\Omega_{k}$ is a homeomorphic image of
${\cal U}$, and is therefore itself a contractible subdomain, hence
$\Omega_{k+1}$ is simply connected. 

\vspace{.1in}

Now consider a simple closed loop $C\subset\Omega_{\infty}$, i.e., 
$C\subset\Omega_{k}$ for all $k$, and hence there exists a contractible
subdomain ${\cal U}_{k}$ in $\Omega_{k}$ such that $\partial{\cal U}_{k}
 = C$. But uniqueness of the interior of $C$ of course implies that 
${\cal U}_{k}\equiv{\cal U}$ for all $k$, and hence ${\cal U}\subset
\Omega_{\infty}$. 

\end{pf}

\begin{definition} The Reeb flow will be said to be 
 ``locally recurrent'' near a periodic orbit ${\cal P}$ if it is 
Lyapunov--stable within a tubular neighbourhood $T_{\cal P}$
and for any sufficiently small 
embedded disc $\tilde{\Delta}$, corresponding to a transversal slice 
through $T_{\cal P}$ at some point $p_{0}$, the 
limit set $\Omega_{\infty}\subseteq\tilde{\Delta'}\subseteq\tilde
{\Delta}$ is open. An
orbit ${\cal P}$ itself may also be referred to as ``locally 
recurrent'' in this context.

\end{definition} 

\vspace{.1in}

As mentioned in section 1, the treatment of a finite--energy 
pseudoholomorphic map $(\psi,a):\Sigma@>>>M\times{\Bbb R}$
from a punctured Riemann surface, asymptotic to a periodic orbit,
may be restricted to a map between $D\setminus\{0\}\subset{\Bbb C}$
and, for convenience, the Martinet neighbourhood $\Delta\times
{\Bbb S}^{1}\times[a_{0},\infty)$ of ${\cal P}$, with coordinates 
$(x,y,\theta,a)$ (cf. [8], [10], [11]). In particular, for $p_{0}$
corresponding to the origin of $\Delta$, we have $\gamma_{0}(t) = 
(0,0,e^{{\bold i}(\theta_{0}+2\pi t/\tau_{0})})$, where $\tau_{0}$ 
denotes the minimal period of ${\cal P}$.
In relation to the asymptotic 
results of [8], it will sometimes be convenient to use cylindrical
coordinates on $D\setminus\{0\}$, viz $(r,\varphi) = (-\ln(|z|),
arg(z))\in [r_{0},\infty)\times{\Bbb S}^{1}$. In particular, the 
{\em charge} of the puncture is given by
\[ T:= \frac{1}{2\pi}\lim_{r@>>>\infty}\int_{{\Bbb S}^{1}}{\psi}^{*}
\lambda \ ,\]
and is a well--defined integer multiple of $\tau_{0}$ (cf. [8]).
Within this coordinate system, the asymptotic relations of Hofer,
Wysocki and Zehnder (mentioned in section one) may be interpreted as
\[a(r,\varphi) = Tr + a_{0} + \varepsilon(r,\varphi) \ ,\]
\[\theta(r,\varphi) = T(\frac{\varphi}{\tau_{0}}) + \theta_{0} +
\delta(r,\varphi) \ ,\]
with $\varepsilon \ , \ \delta$ approaching zero uniformly in $\varphi$
as $r@>>>\infty$. For $|z|$ sufficiently small (i.e., $r$ sufficiently
large) we may think of $\theta^{-1}(\theta_{0})$ as a union of radial
arcs meeting at the origin and differing by only a small perturbation 
from the rays defined by $arg(z) = \frac{1}{n}(\theta_{0} + 2\pi k) \ , 
\ 0\leq k\leq n-1$, where $n = \frac{T}{\tau_{0}}$. Hence we will 
consider $D\setminus\{0\}$ as a union of $n$ {\em quasi--sectors} 
$Q_{k}$, bounded by these arcs, on which the quasi--analytic ``branches'' 
of $(\psi,a)$ are defined. 

\vspace{.1in}

For the local model in section 2 it was sufficient to assume that $\psi$
is a local embedding in order to represent the image as a graph on which
the restriction of an almost complex structure $J$ depends on just two
of the coordinates of a Darboux chart. Under the assumption that the 
Reeb flow is locally recurrent near ${\cal P}$, we now select 
$\Omega_{\infty}\times\{\theta_{0}\}$
as coordinate disc within the initial Martinet tube (on which $\lambda
 = f\cdot(d\vartheta +xdy)$, for a function $f$ such that $f(0,0,\vartheta)
\equiv\tau_{0}$ and $\nabla f(0,0,\vartheta)\equiv{\bold 0}$ as described
in [8],[11], though these facts are not used here). Without loss of 
generality, let $\theta_{0}$ be zero and 
consider the cylinder $\Omega_{\infty}\times[0,2\pi]$, which maps to the 
tube via the obvious identification $mod(2\pi)$. The cylinder has ${\cal P}$ 
as its axis, $x=y=0$, and the Reeb vector field in Martinet coordinates 
already looks like $\frac{1}{\tau_{0}}\frac{\partial}{\partial\theta}$ 
when restricted to ${\cal P}$. There is
no consequent loss of generality if we ``normalise'' $\lambda$ by the
constant multiple $\frac{1}{\tau_{0}}$, so that the minimal period is
effectively $1$, and hence the charge $T$ is an integer. By analogy with
the standard construction of Darboux coordinates, the next step is to 
define 
\[{\cal C}:=\{(p,t) \ | \ p=(x,y)\in\Omega_{\infty} \ , \ \hspace{
.1in} 0\leq t\leq\tau(p) \ \} \ ,\]
and a homeomorphism 
\[h:{\cal C}@>>>\Omega_{\infty}\times[0,2\pi] \ , \hspace{.1in}
h\mid_{\Omega_{\infty}\times\{0\}} = 1 \ ,\]
which is smooth for all $0<t<\tau(p)$, coming from solutions of the
ordinary differential equation
\[\frac{d\gamma_{p}}{dt} = X_{\lambda}(\gamma_{p}(t)) \ .\]
It follows that on the interior of ${\cal C}$, the standard contact form
$\lambda_{0}$ and $\lambda':= h^{*}\lambda$ have the same Reeb vector field, 
corresponding to $\frac{\partial}{\partial t}$. 
We now consider the Cauchy--Riemann system
\[\pi((h^{-1}\psi)_{\eta}) + J\pi((h^{-1}\psi)_{\zeta}) = 0 \ 
\hspace{.1in}(**) \ ,\]
\[\lambda'((h^{-1}\psi)_{\zeta}) = -a_{\eta} \ \hspace{.1in} \ 
\lambda'((h^{-1}\psi)_{\eta}) = a_{\zeta} \ \hspace{.1in}
(\dag\dag) \ .\]
As in the local model, for a sufficiently ``thin'' neighbourhood of
${\cal P}$, the standard projection $(v_{1},v_{2},v_{3})\mapsto (v_{1},
v_{2})$ determines a linear isomorphism $\mu$ between $\xi' := 
ker(\lambda')$ and ${\Bbb R}^{2}$. Hence we define a $2\times 2$
matrix--valued function $j(x,y) = \mu\circ J\circ\mu^{-1}$, such that
${\bold x}:=(x,y)$ implies (**) can be written in the form
\[{\bold x}_{\eta}(z) + j{\bold x}_{\zeta}(z) = 0 \ .\]
Let $\alpha$ denote the diffeomorphism of $\Omega_{\infty}
\times\{\theta_{0}\}$ defined by the return map $\alpha(p):= 
\gamma_{p}(\tau(p))$, hence $\alpha(0) = 0$.

\begin{lemma} If ${\cal L}_{X_{\lambda}}J = 0$, then in a neighbourhood
of $0\in\Delta$, the smooth automorphism $\alpha$ is diffeomorphically
equivalent to a rotation.
\end{lemma}

\begin{pf} Letting ${\bold e}_{1}$ denote the vector $(1,0)$, consider 
the system of ordinary differential equations 
\[\frac{d{\bold x}}{ds} = j({\bold x})\cdot{\bold e}_{1} \ ,\]
and define a local diffeomorphism $\varphi:\Delta''@>>>U\subseteq
\Omega_{\infty}$ via existence and uniqueness of solutions, such that
\[(\dot{x}_{1},\dot{x}_{2}) = j({\bold x})\cdot{\bold e}_{1} \ 
\Rightarrow\varphi(x,s) := (x_{1}(s),x_{2}(s)) \ , \]
$\varphi(x,0) = (x,0) \ .$ Hence 
\[\varphi_{*}^{-1}\circ j\circ\varphi_{*} = j_{0} = 
\left(\begin{array}{cc}
0 & -1 \\
1 & 0
\end{array}\right) \ . \]

Let $\Omega_{\infty}'$
denote the simply connected domain inside $U$ which is stabilised by the Reeb flow.
The diffeomorphism $\hat{\alpha}:= \varphi^{-1}\circ\alpha
\circ\varphi$ \  then acts on $\varphi^{-1}(\Omega_{\infty}')\subseteq
\Delta''$ as a restricted automorphism such that $\hat{\alpha}(0) = 0$
under the assumption of local recurrence. The additional assumption 
${\cal L}_{X_{\lambda}}J = 0$ implies that $\alpha^{*}j = j$, hence in
particular $\hat{\alpha}j_{0} = j_{0}\hat{\alpha}$, i.e., $\hat{\alpha}$
is a conformal automorphism. Modulo a conformal transformation
identifying $\varphi^{-1}(\Omega_{\infty}')$ with a disc, $\hat{\alpha}$
is then equivalent to a rotation. 
\end{pf}
 
Now $\tilde{\bold x}_{\zeta}:=\varphi_{*}^{-1}({\bold x}_{\zeta})$ and 
$\tilde{\bold x}_{\eta}:=\varphi_{*}^{-1}({\bold x}_{\eta})$ implies
\[\tilde{\bold x}_{\eta}(\eta,\zeta) + j_{0}\cdot
\tilde{\bold x}_{\zeta}(\eta,\zeta) = {\bold 0} \ (\dag *) .\]
Note that each ``branch'' of
\[\Psi:=(h\circ(\varphi\times 1))^{-1}\psi\]
is defined smoothly in the interior and continuously up to the 
boundaries of a quasi--sector $Q_{k}$
in $D\setminus\{0\}$, but discontinuities arise at points $z_{0}$
lying on the smooth arcs that bound adjacent sectors (in the usual
way ``$\pm$'' will be used to denote opposite sides of the boundary). 
Discontinuities of the transverse 
projection of $\Psi$ are therefore described by the relations
\[\lim_{z@>>>z_{0}^{\pm}}\tilde{{\bold x}}(z):= \tilde{{\bold x}}^{\pm}
(z_{0}) \ \Rightarrow \ \hat{\alpha}(\tilde{{\bold x}}^{-}(z_{0})) = 
\tilde{{\bold x}}^{+}(z_{0}) \ .\]
Hence on each $Q_{k}\subset D\setminus\{0\} \ , \  (\dag *)$ defines a
holomorphic 
function  $w=F_{k}(z)$ which partially describes a branch of $\Psi$, 
such that
\[F_{k}\mid_{\overline{Q}_{k}\cap\overline{Q}_{k+1}}= \hat{\alpha}
\circ F_{k+1}\mid_
{\overline{Q}_{k}\cap\overline{Q}_{k+1}} \ ,\ \hspace{.1in} 0\leq 
k\leq n-1 \ .\]

\vspace{.1in}

\begin{lemma} If the Reeb flow determines a fibration
of a tubular neighbourhood by periodic orbits over $\Omega_{\infty}$,
i.e., $\alpha = 1$, then there is a single holomorphic function 
$F(z)$ on $D$ \ , \ $F(0) = 0$ \ , \  describing the transverse 
projection of $\Psi$.
\end{lemma}

\begin{pf} $\alpha = 1$ implies the existence of a function
$F(z)$, continuous on $D$ and holomorphic on the interior
of each quasi--sector. The demonstration that $F$ is holomorphic on
$D$ is a standard application of Morera's Theorem (and the Removable 
Singularities Theorem at the origin).
Specifically, let $z_{0}$ lie on one of the smooth arcs bounding a
quasi--sector and let $z' = \eta' + {\bold i}\zeta'$ be a local complex
coordinate with respect to which the arc is locally described as the
graph of a function $\zeta' = \rho(\eta')$. Let $\Gamma$ be a simple
loop inside the $z'$--coordinate neighbourhood of $z_{0}$.
Clearly if $\Gamma$ does not intersect the arc, then $F$ must be 
holomorphic on a slightly larger domain containing $\Gamma$, hence
\[\int_{\Gamma}F \ dz' = 0 \ .\]
If $\Gamma$ intersects the arc, consider a simple affine coordinate
transformation so that the cord joining the endpoints $a,b$ of the 
intersection is now the axis corresponding to $\zeta' = 0$. Hence,
without loss of generality, we may assume $\rho(a) = \rho(b) = 0$.
Let $\gamma$ denote the arc corresponding to the graph between
$a$ and $b$, and for $\varepsilon > 0$ let $\gamma_{\pm\varepsilon}$
denote the arcs corresponding to $\zeta' = \rho(\eta')\pm\varepsilon$,
lying between $a\pm{\bold i}\varepsilon$ and $b\pm{\bold i}\varepsilon$
and on either side of $\gamma$. From the continuity of $F$ it follows
that
\[\int_{\gamma}F \ dz' = \int_{a}^{b}F(\eta' + {\bold i}\rho(\eta'))
d\eta' = \lim_{\varepsilon @>>> 0}\int_{a}^{b}F(\eta' + {\bold i}
(\rho(\eta') \pm \varepsilon)) \ d\eta'\]
\[\hspace{.4in} = \lim_{\varepsilon @>>> 0}
\int_{\gamma_{\varepsilon}}F \ dz' \ .\]
Now decomposing $\Gamma$ into two simple loops $\Gamma_{\pm}$
having $\gamma$ as their common boundary component, and applying the 
above limit, we see that 
\[\int_{\Gamma}F \ dz' = 0 \ .\]
Hence $F$ extends holomorphically to a neighbourhood of $z_{0}$.
\end{pf}
 
Returning now to the particular form of the equations $(\dag\dag)$ in
${\cal C}$ note that any $\lambda'$ with $\frac{\partial}{\partial t}$
as its Reeb vector field must take the general form 
\[\lambda' = dt + f_{1}({\bold x})dx + f_{2}({\bold x})dy \ .\]
As was seen with respect to a Darboux chart of the local model,
such a local presentation of the contact form allows decoupling of
$(\dag\dag)$ into an inhomogeneous Cauchy--Riemann equation. This 
property is preserved under the diffeomorphism $\varphi$, however,
if it is assumed that ${\cal L}_{X_{\lambda}}
J = 0$, hence in particular the matrix $j$ above is independent of $t$.   
In this case, letting $u = t+{\bold i}a$, $(\dag\dag)$ becomes

\[\frac{\partial}{\partial\bar{z}}[u\mid_{Q_{k}}] = -\{f_{1}(F_{k}(z))
\frac{\partial\Re(F_{k})}{\partial\bar{z}} + f_{2}(F_{k}(z))
\frac{\partial\Im(F_{k})}{\partial\bar{z}}\}\]
\[ = -\frac{1}{2}[(f_{1} +{\bold i}f_{2})
\circ F_{k}(z)]\cdot\overline{F'_{k}(z)} \ , \ \]
keeping in mind that this equation is defined smoothly
only on the interior of each quasi--sector $Q_{k}$. Define
$\omega := \frac{1}{2}(f_{1}+{\bold i}f_{2})(w,\bar{w})d\bar{w}$, so that
\[\lambda' = dt + 2\Re(\omega) \ , \ \hbox{and}\]
\[\frac{1}{2}[(f_{1} +{\bold i}f_{2})\circ F_{k}(z)]\cdot\overline{F'_{k}(z)}
d\bar{z} = F_{k}^{*}\omega \ .\]
Now ${\cal L}_{X_{\lambda}}\lambda = 0$ implies $\hat{\alpha}^{*}
\Re(\omega) = \Re(\omega)$. In particular, ${\bf f}:= (f_{1},f_{2})$
implies $\Re(\omega) = ({\bf f},*)$ with respect to the standard inner
product on ${\Bbb R}^{2}$, and hence $\hat{\alpha}_{*}^{t}{\bf f} = 
{\bf f}$. Similarly $\Im(\omega) = (j_{0}{\bf f},*)$, while
${\cal L}_{X_{\lambda}}J = 0$ implies $\hat{\alpha}_{*}^{t}j_{0} = 
j_{0}\hat{\alpha}_{*}^{t}$, so that
\[\hat{\alpha}^{*}\Im(\omega) = 
(\hat{\alpha}_{*}^{t}j_{0}{\bf f},*) = (j_{0}{\bf f},*) = \Im(\omega) \ .\]
It follows that $\hat{\alpha}^{*}\omega = \omega$, and hence
\[F_{k}^{*}\omega\mid_{\overline{Q}_{k}\cap\overline{Q}_{k+1}} = 
F_{k}^{*}(\hat{\alpha}^{*}\omega)\mid_{\overline{Q}_{k}\cap\overline
{Q}_{k+1}} = 
(\hat{\alpha}\circ F_{k})^{*}\omega\mid_{\overline{Q}_{k}\cap
\overline{Q}_{k+1}} \] 
\[ = F_{k+1}^{*}\omega\mid_{\overline{Q}_{k}\cap\overline{Q}_{k+1}} \ .\]
There now exists a 
continuous function $G(z,\bar{z})$ on $D$ such that
\[G(z,\bar{z})d\bar{z}\mid_{Q_{k}}:= F_{k}^{*}\omega \ , \ \hspace{.1in}
0\leq k\leq n-1 \ . \]
Letting $\Phi(w,\bar{w}) = f_{1}+{\bold i}f_{2}$ and assuming, from lemma 1,
that $\hat{\alpha}$ is a rotation, we may write
\[\hat{\alpha}^{*}\omega = \omega \ \Rightarrow \Phi\circ\hat{\alpha}
(w,\bar{w})
 = \overline{\hat{\alpha}^{-1}}\Phi(w,\bar{w}) \ .\]
It is easily seen that
\[\Phi_{w}\circ\hat{\alpha}(w,\bar{w}) = |\hat{\alpha}|^{-2}\Phi_{w} \ , 
\hspace{.1in} 
\Phi_{\bar{w}}\circ\hat{\alpha}(w,\bar{w}) = \overline{\hat{\alpha}^{-2}}
\Phi_{\bar{w}} \ ,\]
and hence 
\[\frac{\partial G}{\partial z}\mid_{\overline{Q}_{k}\cap\overline{Q}_{k+1}}
 = \Phi_{w}(F_{k+1})\cdot |F_{k+1}'|^{2}\]
\[ = \Phi_{w}(\hat{\alpha}\cdot F_{k})
\cdot|\hat{\alpha}\cdot F_{k}'|^{2} = \Phi_{w}(F_{k})\cdot|F_{k}'|^{2} \ ,\]
while
\[\frac{\partial G}{\partial \bar{z}}\mid_{\overline{Q}_{k}\cap\overline
{Q}_{k+1}}
 = \Phi_{\bar{w}}(F_{k+1})\cdot(\bar{F}_{k+1}')^{2} + \Phi(F_{k+1})\cdot
\bar{F}_{k+1}'' \]
\[= \Phi_{\bar{w}}(\hat{\alpha}\cdot F_{k})
\cdot(\overline{\hat{\alpha}\cdot F}_{k}')^{2} + \Phi(\hat{\alpha}\cdot 
F_{k})\cdot\overline{\hat{\alpha}\cdot F}_{k}''\]
\[ = \Phi_{\bar{w}}(F_{k})\cdot
(\bar{F}_{k}')^{2} + \Phi(F_{k})\cdot\bar{F}_{k}'' \ .\]
Thus $G$ is a continuously differentiable function, and 
\[\hat{G}(z,\bar{z}) := \int_{D}\frac{G(\mu,\bar{\mu})}{\mu-z}
d\mu\wedge d\bar{\mu}\]
is a twice-continuously differentiable function on $D$. 
We may therefore write the solutions
$u = H_{k}(z) - \frac{1}{2\pi{\bold i}}\hat{G}(z,\bar{z})$
for a holomorphic function $H_{k}$ defined on each $Q_{k}$. Now
$a(z) = \Im(u)$
is a smooth function on $D\setminus\{0\}$, while $\hat{G}$  belongs
 to $C^{2}(D)$, hence $\Im(H_{k}(z))$ 
must in fact correspond to a single function $h(z)$ for all $k$,
which belongs to $C^{2}(D\setminus\{0\})$ and is harmonic 
inside each $Q_{k}$, therefore harmonic throughout $D\setminus\{0\}$. 
Moreover, the
harmonic conjugate of $h$ is uniquely defined up to a constant, hence
we have a single harmonic function $\hat{h}$ on the 
punctured disc such that 
\[\Re(H_{k}(z)) = \hat{h}(z) + c_{k} \ , \ \hspace{.1in} 0\leq k\leq
n-1 \ .\]

\vspace{.1in}

Recall that $\Re(u)$ corresponds to $t$ such that $0<t<\tau(p)$ for some
$(p,t)$ belonging to the image of $\Psi$. Discontinuities of $t$ along
the boundaries of each $Q_{k}$ are consequently determined by $\tau(p)$,
i.e.,
\[\lim_{z@>>>z_{0}^{-}}\hat{h}_{k}(z) = \lim_{z@>>>z_{0}^{+}}
\hat{h}_{k+1}(z) + \tau(\tilde{{\bold x}}^{-}(z_{0}))\]
for $z_{0}$ lying in the boundary arc $\overline{Q}_{k}\cap\overline{
Q}_{k+1}$.
From the discussion above it follows that 
$\tau(\tilde{{\bold x}}^{-}(z))$ must be constant for $\Psi$ 
restricted to a boundary arc. By
continuity of these arcs as they radiate from the origin of $D$, this 
constant value must correspond to $\tau_{0}$, or $1$
if $\lambda$ is assumed to have been normalised. The arcs themselves
were originally defined by the arbitrary choice of $\theta_{0}\in{\Bbb
S}^{1}$, hence we may conclude that $\tau(p) = 1$ for all $p\in\Omega_
{\infty}'$. Moreover, we have
\[c_{k+1} - c_{k} = 1 \ , \ \hspace{.1in} 0\leq k\leq n-1 \ ,\]
where, without loss of generality, we may set $c_{0} = 0$, hence
$c_{k} = k$.
Recalling the asymptotic formulae $a(z) = -\ln|z|^{n} + \varepsilon(z)$ and
 $\theta(z) = arg(z^{n}) + \delta(z)$, one may notice also the 
approximation of the holomorphic part of $u$ by the branched analytic function 
$log(z^{n})$. Specifically,
\[ u = \frac{1}{2\pi{\bold i}}\{\log(\rho(z)e^{-\hat{G}(z,\bar{z})})\} \ ,\]
where the analytic function  $\rho(z):= e^{2\pi H(z)}$ has 
order $n$ at $z = 0$, and $K:=sup_{D}|G(z,\bar{z})|$ implies
\[\frac{1}{2\pi}|\hat{G}(z,\bar{z})| = \frac{1}{2\pi}|\int_{D}
\frac{G(\mu,\bar{\mu})}{\mu-z}d\mu\wedge d\bar{\mu}|\leq K|z|\]
is uniformly bounded in the parameter $arg(z)$, i.e., 
\[\lim_{r@>>>\infty}|\hat{G}(r,\varphi)| = 0\hspace{.1in}\hbox{in}
\hspace{.1in} C^{0}({\Bbb S}^{1}) \ .\]   
 In summary, we have the following

\begin{theorem} Let $(\psi,a):D\setminus\{0\}@>>>M\times{\Bbb R}$ be
a $J$--holomorphic curve of finite energy and charge $n$ at $z=0$,
asymptotic to a locally recurrent periodic orbit ${\cal P}$, near which 
${\cal L}_{X_{\lambda}}J = 0$. Consider any tubular neighbourhood
of ${\cal P}$ in $M$, diffeomorphic to $\Delta\times{\Bbb S}^{1}$ such
that $\{0\}\times{\Bbb S}^{1}\approx{\cal P}$. There exists a diffeomorphic
change of coordinates in $\Delta\times[0,2\pi)$ such that on each 
quasi--sector $Q_{k}\subset D\setminus\{0\}$ the map $(\psi,a)$ 
can be expressed in the form 
\[(F_{k}(z), H_{k}(z) - \frac{1}{2\pi{\bold i}}\hat{G}(z,\bar{z}))\ , \ 
\hspace{.1in} 0\leq k\leq n-1 ,\]
where $F_{k} \ , \ H_{k}$ are holomorphic on $Q_{k}$ and continuous
on $\overline{Q}_{k}$, such that
\[F_{k}\mid_{\overline{Q}_{k}\cap\overline{Q}_{k+1}} = \hat{\alpha}
\circ F_{k+1}\mid_{\overline{Q}_{k}\cap\overline{Q}_{k+1}} \ , \]
while each $H_{k}$ corresponds to an
analytic branch of $\frac{1}{2\pi{\bold i}}\log(\rho) \ , \ 
ord_{0}(\rho) = n$.
Moreover, the function $\hat{G}$ belongs to $C^{2}(D)$ and is bounded 
by $K|z|$. Finally, if $\alpha = 1$, then each $F_{k}$ is 
the restriction of a single function $F$ holomorphic on $D$, \ $F(0) = 0$.

\end{theorem}

A representation corresponding to the classical local parametrization
of plane algebroid curves is easily obtained as follows for the case 
$\alpha = 1$. Let
\[g(w,\bar{w}) := \frac{1}{2\pi{\bold i}}\int_{\Delta'}\frac{(f_{1}
+{\bold i}f_{2})(\mu,\bar{\mu})d\mu\wedge d\bar{\mu}}{\mu - w} \ ,\]
so that $\bar{\partial}g = \omega$. Then 
\[G(z,\bar{z})d\bar{z} = F^{*}\omega = F^{*}(\bar{\partial}g)\]
\[ = \bar{\partial}(F^{*}g) = \bar{\partial}(g\circ F) \ ,\]
and hence $\hat{G} = 2\pi{\bold i}(g\circ F + \hat{H})$,
for some holomorphic function $\hat{H}$. Now define a coordinate
\[ v:= e^{2\pi{\bold i}u} = \rho(z)e^{-\hat{G}} = \rho(z)e^{-g\circ F 
- \hat{H}}\]
and let $\rho(z)e^{-\hat{H}} = f(z)^{n}$ for some holomorphic function $f$
on $D'\subseteq D$, with $f'(0)\neq 0$. Now $\xi:=f(z)$ implies $w = F(z)
 = F\circ f^{-1}(\xi)$, while $v = \xi^{n}e^{-g\circ F\circ f^{-1}(\xi)}$.
From the final coordinate diffeomorphism $w' = w \ ; \ v' = ve^{g(w,
\bar{w})}$, it now follows that
\[ w' = F\circ f^{-1}(\xi) \ ; \ \hspace{.1in} v' = \xi^{n} \ .\]

\section{$\alpha = 1$: Circle fibrations of ${\Bbb S}^{3}$} 

A class of examples of tight contact structures for which $\alpha = 1$
and ${\cal L}_{X_{\lambda}}J = 0$ near a periodic orbit is provided by
the circle fibrations of ${\Bbb S}^{3}$, of which the most elementary
instances are the Hopf fibrations. Given a four--dimensional real vector 
space $V$, these fibrations correspond to families of invariant planes
(i.e., ``complex lines'') distinguished by linear endomorphisms $J_{0}$
determining standard complex structures on $V$, and are parametrised
by ${\Bbb O}(4)/{\Bbb U}(2)$. The base space of each such fibration is
a (Riemann) sphere inside the Grassmann manifold ${\Bbb G}r_{2}(V)$.
More general fibrations correspond to families of planes (``line
congruences'') parametrised by compact surfaces $\Sigma\subset{\Bbb G}
r_{2}(V)$. A line congruence is said to be {\em elliptic} if for all
$P\in\Sigma$, there exists a 2--sphere ${\Bbb S}$ corresponding to some
$J_{0}$ such that $T_{P}\Sigma = T_{P}{\Bbb S}$ inside ${\Bbb G}r_{2}(V)$,
hence in particular $\Sigma$ is itself diffeomorphic to a sphere.        
Let $\bigwedge^{2} V$ denote the space of exterior 2-forms, on which the
duality operator acts in the standard way. The spaces of ``self-dual''
and ``anti--self--dual'' forms then correspond to +1 and --1 eigenspaces
of this operator, defining a direct sum decomposition 
$\bigwedge^{2} V\cong\bigwedge^{2}_{+}\bigoplus\bigwedge^{2}_{-} \ .$
If ${\Bbb S}_{+}$ and ${\Bbb S}_{-}$ denote the 2--spheres of radius 
$\frac{1}{\sqrt{2}}$ inside each of these eigenspaces, then the Grassmann
manifold of oriented 2--planes of $V$ is well--known to correspond to
${\Bbb S}_{+}\times{\Bbb S}_{-}$. Moreover, it was shown by Gluck and 
Warner [4] that the surfaces $\Sigma$ of generalised Hopf fibrations
are precisely the graphs of distance--decreasing smooth maps
$f:{\Bbb S}_{-}@>>>{\Bbb S}_{+}$, with standard Hopf fibrations 
coresponding to constant maps. In [13], McKay also observed that they 
correspond to general elliptic line congruences and hence determine 
non--linear complex 
structures on $V$ which ``osculate'' with linear structures at each
$P\in\Sigma$. We apply this idea to the explicit construction of 
contact structures on ${\Bbb S}^{3}$ as follows. Let ${\cal J}(V)$
denote the space of linear endomorphisms of $V$ corresponding to 
linear complex structures, then it was shown in [13] that each
elliptic line congruence $\Sigma$ determines a map 
\[J:\Sigma @>>> {\cal J}(V)\subset V\otimes V^{*}\]
such that for each plane $P\in\Sigma$, $J(P)$ is linear and is the 
``osculating'' structure to $\Sigma$ at $P$ in the sense that
both $P$ and $P^{\perp}$ (with respect to a given inner product on $V$)
are complex lines relative to $J(P)$. The family of planes 
determined by $\Sigma$ describes a rank--two vector bundle ${\cal P}
@>\pi>>\Sigma$ such that the total space, corresponding
to the incidence manifold ${\cal P} = \{({\bold v},P)\in V\times\Sigma 
\ | \ {\bold v}\in P\}$ also maps surjectively to $V$. In fact, there 
exists ${\cal P}@>\sigma>>V$ such that ${\cal P}\setminus\sigma^{-1}(0)
\cong V\setminus\{0\}$. Hence define $\varphi:=\pi\circ\sigma^{-1}:V\setminus
\{0\}@>>>\Sigma$, so that $\varphi^{-1}(P) = \{{\bold v}\in V\setminus
\{0\} \ | \ {\bold v}\in P\}$ \ . 

Consider the pullback $\varphi^{*}J:V\setminus\{0\}@>>>V\otimes V^{*}$.
Hence with respect to a designated orthonormal basis of $V$, noting that
$J$ and $\varphi^{*}J$ are skew--symmetric matrix--valued functions, we
may represent it in the form
\[\varphi^{*}J =\Sigma_{\mu.\nu}J^{\mu\nu}\frac{\partial}{\partial x
_{\mu}}\wedge dx_{\nu} \ .\]
Note moreover that the isomorphism $V^{*}\cong V$ via the Euclidean 
inner product allows us to 
define a 2--form
$\omega_{\Sigma}:= \Sigma_{\mu.\nu}J^{\mu\nu}dx_{\mu}\wedge dx_{\nu}$.
In the following, let $i_{\bold v}$ denote contraction of a form by the
position vector ${\bold v} = \frac{1}{2}\Sigma_{\eta}v_{\eta}\frac{
\partial}{\partial\eta}$ \ .

\begin{lemma} If $d\omega_{\Sigma} = 0$ then $\omega_{\Sigma}({\bold v})
 = d(i_{\bold v}\omega_{\Sigma})$ \ .
\end{lemma}

\begin{pf} 
\[i_{\bold v}\omega_{\Sigma} = \frac{1}{2}\Sigma_{\mu < \nu}J^{\mu\nu}
(x_{\mu}dx_{\nu} - x_{\nu}dx_{\mu}) \]
\[ = \frac{1}{2}\Sigma_{\beta\neq\alpha}(-1)^{\varepsilon}x_{\beta}
J^{\alpha\beta}dx_{\alpha} \ ,\]
where 
$$\varepsilon = 
\left\{\begin{array}{l} 0 \hspace{.1in} \beta < \alpha \\
 1 \hspace{.1in} \beta > \alpha \ . \\
\end{array} \right. \\ $$

Therefore
\[d(i_{\bold v}\omega_{\Sigma}) = \omega_{\Sigma} + \frac{1}{2}\Sigma_
{\gamma < \alpha\neq\beta}(-1)^{\varepsilon}x_{\beta}(\frac{\partial
J^{\alpha\beta}}{\partial x_{\gamma}} - \frac{\partial
J^{\gamma\beta}}{\partial x_{\alpha}})dx_{\gamma}\wedge dx_{\alpha} \]
\[ = \omega_{\Sigma} + \frac{1}{2}\Sigma_{\gamma < \alpha\neq\beta}
x_{\beta}\frac{\partial
J^{\gamma\alpha}}{\partial x_{\beta}}dx_{\gamma}\wedge dx_{\alpha} \]
(using the relations provided by $d\omega_{\Sigma} = 0$)
\[ = \omega_{\Sigma} + \frac{1}{2}\Sigma_{\gamma < \alpha}
(\nabla_{\bold v}J^{\gamma\alpha})dx_{\gamma}\wedge dx_{\alpha} \ .\]
Note that the functions $J^{\gamma\alpha}$, obtained by pulling back 
$J$, are constant in the radial directions of $V$, hence the directional
derivatives $\nabla_{\bold v}J^{\gamma\alpha} = 0$ for all ${\bold v}
\in V\setminus\{0\}$ \ . 

\end{pf}
     
Now define the 1--form $\lambda_{\Sigma}:= i_{\bold v}\omega_{\Sigma}$
for $\omega_{\Sigma}$ closed,
noting that if $[J]$ denotes the matrix of $\varphi^{*}J$, then we may 
express $\lambda_{\Sigma}$ in terms of the inner product as  
\[\lambda_{\Sigma}({\bold w}) = - ({\bold w},[J]\cdot{\bold v}) \]
for all ${\bold w}\in V$. Moreover, let $X_{\lambda}$ be the vector 
field defined by $X_{\lambda}({\bold v}) = -[J]\cdot{\bold v}$, so that

$(i) \ X_{\lambda}$ is tangent to ${\Bbb S}^{3}$, since
\[ ({\bold v}, X_{\lambda}({\bold v})) = -({\bold v},[J]\cdot{\bold v})
 = -([J]^{t}\cdot{\bold v},{\bold v}) = ([J]\cdot{\bold v},{\bold v})
 \ ,\]
hence $({\bold v}, X_{\lambda}({\bold v})) = 0$ \ ,

$(ii) \ \lambda_{\Sigma}(X_{\lambda}) = |[J]\cdot{\bold v}|^{2} = 1$ \ ,
since $[J]_{\bold v}\in{\Bbb O}(4)$ for each ${\bold v}\in 
V\setminus\{0\}$, and

$(iii)$ for all ${\bold w} \ , \ {\bold u}\in V$ we have
\[\omega_{\Sigma}({\bold w},{\bold u}) = \Sigma_{\mu < \nu} J^{\mu\nu}
(w_{\mu}u_{\nu} - w_{\nu}u_{\mu}) = ({\bold w},[J]\cdot{\bold u}) \ ,\]
and hence
\[i_{X_{\lambda}}\omega_{\Sigma} = -(*,[J]^{2}\cdot{\bold v}) = (*,{\bold v}) 
= 0 \ \]
when restricted to $T{\Bbb S}^{3}$. Clearly $\varphi^{*}J$ is 
preserved by the Reeb flow, i.e., ${\cal L}_{X_{\lambda}}\varphi^{*}J 
= 0$. Moreover, if $P_{\bold v}$ denotes
the subspace spanned by $\{{\bold v} \ , \ [J]\cdot{\bold v}\}$, then
$ker(\lambda)\cap T{\Bbb S}^{3} = P_{\bold v}^{\perp}$ is also an invariant
subspace of $[J]_{\bold v}$, so we may write $j:= [J]\mid_{P_{\bold v}^
{\perp}}$ for all ${\bold v}\in{\Bbb S}^{3}$. In particular, for all 
$\xi\in ker(\lambda)$, we have $\omega_{\Sigma}(\xi,j\cdot\xi) = 
-|\xi|^{2}$, which means that $\lambda_{\Sigma}$ is a contact structure
compatible with the partial complex structure $j$. 

It should be mentioned that all Hopf fibrations, including the non--linear
ones, are smoothly equivalent as circle bundles (cf. [4]), and yet at the
level of contact structures they are distinct. For suppose $\delta:
{\Bbb S}^{3}@>>>{\Bbb S}^{3}$ is a diffeomorphism that identifies the Reeb
flows of a given structure $\lambda_{\Sigma}$ and that of the standard 
structure $\lambda_{0}$. In particular, suppose that $\delta_{*}
X_{\lambda} = X_{\lambda_{0}}$, and $\delta_{*}\circ \varphi^{*}J
\circ\delta_{*}^{-1} = J_{0}$, where $J$ denotes the osculating complex 
structure associated
with $\lambda_{\Sigma}$. Now assume in addition that $\delta^{*}\lambda
_{0} = \lambda_{\Sigma}$, hence
\[\delta^{*}\lambda_{0}({\bold u}) = \lambda_{0}(\delta_{*}{\bold u})
  = (\delta_{*}{\bold u} \ , \ J_{0}\cdot{\bold v}) = \lambda_{\Sigma}
({\bold u}) = ({\bold u} \ , \ \varphi^{*}J\cdot\delta^{-1}_{*}
{\bold v}) \ ,\]
therefore
\[(\delta_{*}{\bold u},J_{0}\cdot{\bold v}) = ({\bold u},\delta_{*}^{-1}
J_{0}\delta_{*}\delta_{*}^{-1}{\bold v}) = ({\bold u},\delta_{*}^{-1}
J_{0}\cdot{\bold v}) \ , \]
that is,
\[({\bold u},\delta_{*}^{t}J_{0}\cdot{\bold v}) = ({\bold u},\delta_{*}
^{-1}J_{0}\cdot{\bold v}) \ ,\]
and hence $\delta_{*}^{t}X_{\lambda_{0}} = \delta_{*}^{-1}X_{\lambda_{0}}$
(or, conversely, $\delta_{*}X_{\lambda} = (\delta_{*}^{t})^{-1}
X_{\lambda}$) \ . On the other hand, $(\delta_{*}{\bold u}, J_{0}\cdot
{\bold v}) = ((\delta_{*}^{-1})^{t}{\bold u}, J_{0}\cdot{\bold v})$
if and only if 
\[(\delta_{*}{\bold u}, \delta_{*}\varphi^{*}J\delta_{*}^{-1}{\bold v}) = 
((\delta_{*}^{-1})^{t}{\bold u}, \delta_{*}\varphi^{*}J\delta_{*}^{-1}
{\bold v})\]
i.e., letting ${\bold v}':= \delta_{*}^{-1}{\bold v}$ \ ,
\[(\delta_{*}^{t}\delta_{*}{\bold u}, \varphi^{*}J\cdot{\bold v}') = 
({\bold u}, \varphi^{*}J\cdot{\bold v}') \ .\]
In particular, we see that $ker(\lambda_{\Sigma})$ is an invariant
subspace of $\delta_{*}^{t}\delta_{*}$. Note moreover that
\[\delta_{*}^{t}\delta_{*}\varphi^{*}J\delta_{*}^{-1}(\delta_{*}^{t})^{-1}
 = (\delta_{*}^{-1}J_{0}^{t}\delta_{*})^{t} = \varphi^{*}J \ ,\]
which implies that $\delta_{*}^{t}\delta_{*}\mid_{ker(\lambda)}$
is complex--linear. But since it is clearly symmetric, it follows that
$\delta_{*}^{t}\delta_{*}\mid_{ker(\lambda)} = c\cdot I$ for some
$c\in{\Bbb R}\setminus\{0\}$. Now $d(\delta^{*}\lambda_{0}) = 
\delta^{*}(d\lambda_{0}) = d\lambda_{\Sigma}$ implies
\[d\lambda_{\Sigma}({\bold u},{\bold w}) = d\lambda_{0}(\delta_{*}{\bold u},
\delta_{*}{\bold w}) = (\delta_{*}{\bold u}, J_{0}\cdot\delta_{*}
{\bold w}) \]
\[ = ({\bold u}, \delta_{*}^{t}J_{0}\delta_{*}{\bold w}) = 
({\bold u}, \delta_{*}^{t}\delta_{*}\varphi^{*}J\cdot{\bold w}) = 
({\bold u}, \varphi^{*}J\delta_{*}^{t}\delta_{*}{\bold w}) = 
c({\bold u},\varphi^{*}J\cdot{\bold w})\]
(if ${\bold w}\in ker(\lambda_{\Sigma})$), and hence
\[d\lambda_{\Sigma}({\bold u},{\bold w}) = c\cdot d\lambda_{\Sigma}
({\bold u}, {\bold w}) \ , \hspace{.1in} \hbox{i.e.} \ c = 1 \ .\]
We conclude that $\delta_{*}^{t} = \delta_{*}^{-1}$, and hence that
$\delta\in{\Bbb O}(4)$, which restricts any such equivalence of contact
structures to the family of linear Hopf fibrations. In summary:

\begin{theorem} Let $J$ be the osculating complex structure of an
elliptic line congruence $\Sigma\subset{\Bbb G}r_{2}(V)$, such that
the skew--symmetric 2--form $\omega_{\Sigma}$ on $V\setminus\{0\}$
is closed. Then the 1--form $\lambda$, such that $\lambda_{{\bold v}}
:=i_{{\bold v}}\omega_{\Sigma}$, defines a fillable, hence tight,
contact structure on ${\Bbb S}^{3}$ for which ${\cal L}_{X_{\lambda}}J
 = 0$ and $\alpha = 1$. Moreover, two such structures are equivalent
via a diffeomorphism $\delta$ of ${\Bbb S}^{3}$ if and only if 
$\delta\in{\Bbb O}(4)$.

\end{theorem}
 
\section{The classical models of the sphere and rational ellipsoids}

As an explicit illustration of Theorem 1, we will first examine the
classical case of analytic curves in ${\Bbb C}^{2}$. The realization
of algebraic curves as finite-energy pseudoholomorphic maps has
also been examined in [7], though from a slightly different point
of view. Let $M = {\Bbb S}^{3}\subseteq{\Bbb C}^{2}$, with 
\[(\lambda_{0})_{{\bold v}} = (*,J_{0}\cdot{\bold v})\mid_{T{\Bbb S}^{3}}\]
the standard contact form defined with respect to the complex structure
$J_{0}$ of ${\Bbb C}^{2}$ and the Euclidean inner product of ${\Bbb R}^{4}$.
$M$ is then defined by the equation $|w_{1}|^{2}+|w_{2}|^{2} = 1$ with
respect to complex coordinates in ${\Bbb C}^{2}$, while
\[\lambda_{0} = \Re(-{\bold i}(\bar{w}_{1}dw_{1} + \bar{w}_{2}dw_{2})) \ .\]
Consider $\Phi(z):D@>>>{\Bbb C}^{2}$ a complex--analytic curve defined 
on a neighbourhood of the origin in ${\Bbb C}$ such that $\Phi(z) = 
(z^{n},\rho(z))$ and $ord_{0}(\rho)\geq n+1$, hence $\Phi$ has a singularity
of multiplicity $n$ at $(0,0)$. The corresponding map $(\psi, a):D@>>>{\Bbb
S}^{3}\times{\Bbb R}$ is given by 
\[\psi(z) = [|\rho|^{2}+|z|^{2n}]^{-\frac{1}{2}}(z^{n},\rho(z)) \ , \ \]
\[a(z) = -\frac{1}{2}\ln(|z|^{2n}+|\rho|^{2}) = -n\ln|z| - \frac{1}{2}
\ln(1+|z|^{-2n}|\rho|^{2}) \ .\]
The periodic orbit ${\cal P}$ corresponds simply to the circle defined
by $\{w_{2} = 0\}\cap{\Bbb S}^{3}$ and is associated with the degenerate
tangent cone of $\Phi$ at $(0,0)$. It is, moreover, a simple calculation
to verify that $charge_{z=0}(\psi) = n$. Consider the proper holomorphic
map $\sigma:{\cal O}_{{\Bbb P}_{1}}(-1)@>>>{\Bbb C}^{2}$ , where ${\cal O}
_{{\Bbb P}_{1}}(-1)$ denotes the complex line bundle of Chern class
equal to --1 on the Riemann Sphere ${\Bbb P}_{1}({\Bbb C}) = \sigma^{-1}({
\bold 0})$, with $\sigma:{\cal O}_{{\Bbb P}_{1}}(-1)\setminus\sigma^{-1}
({\bold 0})\cong{\Bbb C}^{2}\setminus\{0\}$, given in local coordinates
by the quadratic transformation $w_{1} = \mu \ , \ w_{2} = \mu\nu$. Noting
that $|\mu|^{2}(1+|\nu|^{2}) = 1$ on the chart of $\sigma^{-1}
({\Bbb S}^{3})$ corresponding to $\nu\neq\infty$, we have
\[\sigma^{*}\lambda_{0} = \Re(-{\bold i}(\mu^{-1}d\mu + \bar{\nu}|\mu|^{2}
d\nu)) \ ;\]
moreover, for each $\nu\in{\Bbb P}_{1}({\Bbb C})\setminus\{\infty\}$, 
the Hopf fibration corresponding to 
\[\varpi:\sigma^{-1}({\Bbb S}^{3})@>>>{\Bbb P}_{1}({\Bbb C})
\approx{\Bbb S}^{2}\]
 has fibres $\varpi^{-1}(\nu) = \{\mu\in{\Bbb C} \ | \ 
|\mu| = (1+|\nu|^{2})^{-\frac{1}{2}}\}$. Let $\nu = x+{\bold i}y$, then
$\mu = |\mu|e^{{\bold i}\vartheta} = (1+x^{2}+y^{2})^{-\frac{1}{2}}
e^{{\bold i}\vartheta} \ , \ h:=\sigma\mid_{\varpi^{-1}({\Bbb P}_{1}
\setminus\{\infty\})}$ implies
\[h^{*}\lambda_{0} = \Re(-{\bold i}\{\mu^{-1}d\mu + \bar{\nu}|\mu|^{2}
d\nu\})\]
\[ = \Re(-{\bold i}\{(\sqrt{1+x^{2}+y^{2}} \ )e^{{\bold i}\vartheta}
\left(\frac{-(xdx + ydy)e^{{\bold i}\vartheta}}{\sqrt{(1+x^{2}+y^{2})^{3}}} + 
\frac{{\bold i}e^{{\bold i}\vartheta}d\vartheta}{\sqrt{1+x^{2}+y^{2}}}
\right) \]
\[+ (x-{\bold i}y)\frac{dx+{\bold i}dy}{1+x^{2}+y^{2}}\})\]
\[ = \Re(-{\bold i}\left\{\frac{-(xdx+ydy)}{1+x^{2}+y^{2}} + \frac{xdx+ydy}
{1+x^{2}+y^{2}} + {\bold i}(d\vartheta + \frac{xdy-ydx}{1+x^{2}+y^{2}})
\right\})\]
 \[ = d\vartheta + (1+x^{2}+y^{2})^{-1}(xdy - ydx) \ ,\]
while $h^{*}X_{\lambda_{0}} = \frac{\partial}{\partial\vartheta}$. Now 
\[\Psi(z):=h^{-1}\psi(z) = (\arg(z^{n}) \ ,
 \ F(z)) \hspace{.1in}\hbox{where} \ F(z):=z^{-n}\rho(z) \ .\]
Clearly $\nu = F(z)$ is holomorphic, so it remains to show that
\[h^{*}\lambda_{0}(\Psi_{\eta}) = a_{\zeta} \ ; \ h^{*}\lambda_{0}
(\Psi_{\zeta}) = -a_{\eta} \ .\]

Note that $a(z) = -\frac{1}{2}\ln((\eta^{2}+\zeta^{2})^{n}(1+x^{2}+y^{2}))$
implies
\[a_{\zeta} = \frac{-n\zeta(\eta^{2}+\zeta^{2})^{n-1}(1+x^{2}+y^{2}) - 
(\eta^{2}+\zeta^{2})^{n}(xx_{\zeta}+yy_{\zeta})}{(\eta^{2}+\zeta^{2})^{n}
(1+x^{2}+y^{2})}\]
\[ = \frac{-n\zeta}{\eta^{2}+\zeta^{2}} - 
\frac{xx_{\zeta} + yy_{\zeta}}{1+x^{2}+y^{2}} \ ,\]
while $\Psi_{\eta} = (x_{\eta},y_{\eta},\vartheta_{\eta}) = (
x_{\eta},y_{\eta},\frac{-n\zeta}{\eta^{2}+\zeta^{2}})$ implies 
\[h^{*}\lambda_{0}(\Psi_{\eta}) = \frac{-n\zeta}{\eta^{2}+\zeta^{2}} + 
\frac{xy_{\eta} - yx_{\eta}}{1+x^{2}+y^{2}} \ .\]
But $\nu = x+{\bold i}y = F(z)$ is holomorphic, hence $x_{\eta} = y_{\zeta}
 \ , \ x_{\zeta} = -y_{\eta}$ yields $a_{\zeta} = h^{*}\lambda_{0}
(\Psi_{\eta})$, and similarly for $a_{\eta}$. 

Recalling the discussion of section 3, if we write $\lambda' = dt + 
f_{1}(x,y)dx + f_{2}(x,y)dy$, then $h^{*}\lambda_{0} = 2\pi\lambda'$,
where 
\[ f_{1}(x,y) = \frac{-y}{2\pi(1+|\nu|^{2})} \hspace{.1in} ; \hspace{.1in}
f_{2}(x,y) = \frac{x}{2\pi(1+|\nu|^{2})} \ .\]
Now
\[\frac{1}{2}(f_{1} + {\bold i}f_{2}) = \frac{{\bold i}\nu}
{4\pi(1+|\nu|^{2})} = \frac{{\bold i}}{4\pi}\bar{\partial}\ln(1+|\nu|^{2})
 \ ,\]
and hence $g\circ F = \frac{{\bold i}}{2\pi}\ln(1+|F|^{2})^{\frac{1}{2}}$.
Moreover, we can simply define $\hat{G} = 2\pi{\bold i}g\circ F 
 = -\ln(1+|F|^{2})^{\frac{1}{2}}$. 

On the other hand, modulo rescaling by $2\pi$ so that 
\[\hat{a}(z) = \frac{-1}{2\pi}\ln(|z|^{n}(1+|F|^{2})^{\frac{1}{2}}) \ ,\]
we have
\[ u = t+{\bold i}\hat{a} = \frac{1}{2\pi}(arg(z^{n}) + 
{\bold i}\{-\ln|z|^{n} - \ln(1+|F|^{2})^{\frac{1}{2}}\})\]
\[ = \frac{1}{2\pi{\bold i}}(\log(z^{n}) - \hat{G}) \ ,\]
in accordance with the statement of Theorem 1.   

Let $\Gamma\subset{\Bbb C}^{2}$ be the locus of the plane curve parametrized
by $\Phi$. We remark in conclusion that although the multiplicity of the 
singular point of the strict transform $\overline{\sigma^{-1}(\Gamma)\setminus
\sigma^{-1}(0)}$ is less than $n$, 
the charge at $z=0$ of $\Psi$ is easily seen to be conserved by the 
diffeomorphism $h$. The singular plane curve corresponding to the strict 
transform (assuming the singularity has not been resolved by a
single quadratic transformation) is in fact asymptotic (viewed locally
as a $J$--holomorphic curve $\Psi'$) to a distinct periodic orbit within 
a new 3--sphere bounding a neighbourhood of the transformed singularity. 
However, the linking of the transform of $K_{\Gamma}$ with the original 
periodic orbit of $\sigma^{-1}({\Bbb S}^{3})$ is topologically unaffected. 
   
Now let us turn to the ellipsoids
\[{\Bbb E}_{p,q} := \{(w_{1},w_{2})\in{\Bbb C}^{2} \ | \ p|w_{1}|^{2} +
q|w_{2}|^{2} = 1 \ ; \ (p,q)\in{\Bbb R}^{2}_{+}\} \ .\]
The restriction to ${\Bbb E}_{p,q}$ of $\lambda_{0}$ as defined above 
determines a different contact structure on $T{\Bbb E}_{p,q}$, and in
particular a Reeb vector field $X_{P}:= A\cdot J_{0}\cdot{\bold v}_{P}$
for all $P\in{\Bbb E}_{p,q}$, where 
\[A = \left(\begin{array}{cc}
p & 0 \\
0 & q \\
\end{array}\right)\]
is viewed as an element of ${\bold G}{\bold L}(2,{\Bbb C})$. Solutions
of the equation $\dot{\gamma}(t) = X(\gamma(t))  \ , \ \gamma(0) = P
= (z_{1},z_{2})$, then correspond to maps $t\mapsto(z_{1}e^{{\bold i}pt},
z_{2}e^{{\bold i}qt})$ (cf. e.g., [7]). Note that there are two periodic 
orbits 
corresponding to $w_{1} = 0$ and $w_{2} = 0$ separately. These are the
{\em only} periodic orbits of the Reeb flow if $\frac{p}{q}$ is irrational,
whereas ${\Bbb E}_{p,q}$ is foliated by periodic orbits if $\frac{p}{q}\in
{\Bbb Q}$. Moreover the periodic orbit corresponding to $w_{1} = 0$ has
minimal period $\tau = \frac{2\pi}{q}$ and for $w_{2} = 0$ it is
$\frac{2\pi}{p}$, while $\tau = \frac{2\pi k}{p} = \frac{2\pi l}{q}$ for
all other orbits with respect to fixed relatively prime positive integers
$k,l$. Without loss of generality, consider a tubular neighbourhood 
${\cal U}$ of the
orbit $w_{2} = 0$, and let $\nu:=\frac{z_{2}^{k}}{z_{1}^{l}}$ such that
$l > k$. All periodic orbits of the Reeb vector field within this tubular
neighbourhood then correspond to intersections of ${\Bbb E}_{p,q}$ with
algebraic curves (uniquely determined by $\nu$) of the form 
\[w_{2}^{k} = \nu w_{1}^{l} \ , \ \hspace{.1in} |\nu| < \varepsilon \ ,\]
for some positive $\varepsilon$. Once again, we have
\[\lambda_{0} = \Re(-{\bold i}(\bar{w}_{1}dw_{1} + \bar{w}_{2}dw_{2})) \ ,\]
but this time it will be convenient to introduce a formal coordinate 
transformation of the form
\[w_{1} = \mu \ , \ \hspace{.1in} w_{2} = (\nu\mu^{l})^{\frac{1}{k}}\]
so that
\[\bar{w}_{1}dw_{1} + \bar{w}_{2}dw_{2} = \bar{\mu}d\mu + \frac{1}{k}
|\nu|^{\frac{2}{k}}|\mu|^{\frac{2l}{k}}(\frac{1}{\nu}d\nu + \frac{l}{\mu}
d\mu) \ .\]
Moreover, $p|w_{1}|^{2} + q|w_{2}|^{2} = 1$ implies 
\[|\nu|^{\frac{2}{k}}|\mu|^{\frac{2l}{k}} = \frac{1}{q}(1 - p|\mu|^{2}) 
\ ,\]
hence
\[\bar{w}_{1}dw_{1} + \bar{w}_{2}dw_{2} = \left(\bar{\mu} + \frac{(1-p|
\mu|^{2})l}{kq\mu}\right)d\mu + \left(\frac{1-p|\mu|^{2}}{kq\nu}\right)
d\nu \ .\]
Now $\mu = re^{{\bold i}\vartheta}$ implies $d\mu = e^{{\bold i}
\vartheta}(dr + {\bold i}rd\vartheta)$, so that
\[\left(\bar{\mu} + \frac{(1-p|\mu|^{2})l}{kq\mu}\right)d\mu = 
\frac{1}{p}(\frac{dr}{r} + {\bold i}d\vartheta) \ , \]
while $\nu = x + {\bold i}y$ implies
\[\frac{1-p|\mu|^{2}}{kq\nu}d\nu = \frac{1-p|\mu|^{2}}{kq}\left\{\frac{
xdx + ydy}{x^{2}+y^{2}} + {\bold i}\left(\frac{xdy - ydx}{x^{2}+y^{2}}\right)
\right\} \ ,\]
and hence
\[\lambda_{0} = \frac{1}{p}\{d\vartheta + \frac{(1-p|\mu|^{2})(xdy-ydx)}
{l(x^{2}+y^{2})}\} \ (*).\]

Consider $f(r) := pr^2 + q|\nu|^{\frac{2}{k}}r^{\frac{2l}{k}} - 1$, so that
\[f'(r) = 2pr(1 + \frac{q^2}{p^2}|\nu|^{\frac{2}{k}}r^{2(\frac{l}{k}-1)}) 
\ .\]
Hence $f'(r) = 0$ when $r=0$ or $(\frac{-p^2}{q^{2}|\nu|^{\frac{2}{k}}})^{
\frac{k}{2(l-k)}}$ \ . Moreover $f(0) = -1 \ , \ \lim_{r\rightarrow\infty}
f(r) = \infty$ implies that the equation $f(r) = 0$ has a unique positive
real solution. In other words, the value of $|\mu|$ satisfying the equation
$p|\mu|^{2} + q|\nu|^{\frac{2}{k}}|\mu|^{\frac{2l}{k}} = 1$ is uniquely
determined by $p,q,|\nu|$. Hence write $|\mu| = \varphi(p,q,|\nu|)$
and recall that
\[w_{2} = (\nu(\varphi(|\nu|)e^{{\bold i}\vartheta})^{l})^{\frac{1}{k}}
\hspace{.1in} \hbox{i.e.,}\hspace{.1in} w_{2}^{k}e^{-{\bold i}l\vartheta} = 
\nu\varphi(|\nu|)^{l} \ .\]
Note also that $\varphi(p,q,0) = \frac{1}{\sqrt{p}}$ implies that the complex 
function $\chi(\nu):=\nu\varphi(|\nu|)^{l}$ admits a locally differentiable
inverse, and we may write 
\[\nu = \chi^{-1}(w_{2}^{k}e^{-{\bold i}l\vartheta}) \ .\]
Let $\Delta = \{|w_{2}|<\varepsilon\} \ , \ \Delta' = \{|\nu|<
\varepsilon'\}$ and consider the $k$--fold covering map 
$\beta: {\cal U}\approx\Delta\times{\Bbb S}^{1}@>>>\Delta'\times{\Bbb S}
^{1}$ such that
\[\beta(w_{2},\vartheta) = (\chi^{-1}(w_{2}^{k}e^{-{\bold i}l\vartheta}),
\vartheta) \ .\]
Alternatively, the map $\Theta:{\Bbb C}^{2}@>>>{\Bbb C}^{2}$ such that
$(\mu,\nu) = \Theta(w_{1},w_{2}) = (w_{1}, w_{2}^{k}w_{1}^{-l})$ is 
holomorphic away from
$\{w_{1} = 0\}$ such that $\beta = \Theta\mid_{{\Bbb E}_{p,q}}$.
The equation (*) above may then be written more precisely in the form
$\lambda_{0} = \frac{1}{p}\beta^{*}\lambda$, where 
\[\lambda = d\vartheta + \frac{(1-p\varphi(p,q,|\nu|)^{2})(xdy-ydx)}
{l|\nu|^{2}} \]
(note also that $\lim_{|\nu|\rightarrow 0}\frac{1-p\varphi^2}{l|\nu|^2}
(xdy-ydx) = 0$). Consequently finite--energy pseudoholomorphic maps 
\[(\psi,a):D\setminus\{0\}@>>>\Delta\times{\Bbb S}^{1}\times{\Bbb R}\]
project onto finite--energy maps $(\beta\circ\psi,a)$, pseudoholomorphic with 
respect to $\frac{1}{p}\lambda$. This 
claim is easily verified if we note that 
\[ a_{\zeta} = \frac{1}{p}\beta^{*}\lambda(\psi_{\eta}) = \frac{1}{p}
\lambda(\beta_{*}\circ\psi_{\eta}) = \frac{1}{p}\lambda((\beta\circ\psi)
_{\eta}) \ ,\]
and similarly for $a_{\eta}$. Moreover, 
\[\tilde{\pi}((\beta\circ\psi)_{\eta}) := \beta_{*}\psi_{\eta} - \frac{1}
{p}\lambda((\beta\circ\psi)_{\eta})\tilde{X} \ ,\]
where $\tilde{X}:=\beta_{*}X_{\lambda_{0}} = p\frac{\partial}{\partial
\vartheta}$ is a well--defined vector
field under the above conditions. Hence
\[\tilde{\pi}((\beta\circ\psi)_{\eta}) = \beta_{*}(\pi(\psi_{\eta})) \ 
,\] and therefore
\[0 = \beta_{*}(\pi(\psi_{\eta})+J_{0}\pi(\psi_{\zeta})) = \tilde{\pi}(
\beta\circ\psi)_{\eta} + \beta_{*}(J_{0}\pi(\psi)_{\zeta}) \ .\]
Let $J_{0}$ denote the standard complex structure on ${\Bbb C}^{2}$
as represented by both $(w_{1},w_{2})$ and $(\mu,\nu)$ (and its restriction
to the contact planes of $\lambda_{0}$ and $\lambda$ respectively). Then 
$\beta = \Theta\mid_{{\Bbb E}_{p,q}}$ for $\Theta$ holomorphic implies 
$\beta_{*}\circ J_{0} = J_{0}\circ\beta_{*}$, and hence
\[0 = \tilde{\pi}(\beta\circ\psi)_{\eta} + J_{0}\tilde{\pi}(\beta\circ
\psi)_{\zeta} \ .\]
Note that the tubular neighbourhood into which $\beta\circ\psi$ maps
is fibred by Reeb orbits, hence the return map $\alpha = 1$. By 
comparison, the Reeb flow in a neighbourhood of the original orbit $\{w_{2}
 = 0\}$ in ${\Bbb E}_{p,q}$ induces a return map that is equivalent to a
rational rotation through $2\pi\frac{l}{k}$. In order to find a class of
$J$--holomorphic curves asymptotic to the given periodic orbit within ${\Bbb
E}_{p,q}$, let us first make a harmless renormalization of the contact
structure, i.e, $\lambda_{0}' := p\cdot\lambda_{0} \ , \ X_{\lambda_{0}'}
 = \frac{1}{p}X_{\lambda_{0}}$. Now
\[\lambda_{0}' = d\vartheta + \frac{(1-p\cdot\varphi(|\nu|)
^{2})(xdy-ydx)}{l\cdot|\nu|^{2}}\]
and we are ready to work backwards from a holomorphic parametrization of
the form
\[ \mu = z^{n} \ ; \hspace{.1in} \nu = \Phi(z) \ ,\]
such that $ord_{0}(\Phi) = b\cdot l\geq -nl+1$ and $n = c\cdot k$
 for some integers $b,c$. Now
\[w_{1} = z^{n} \ ;\hspace{.1in} w_{2} = (z^{nl}\cdot\Phi(z))^{\frac{1}{k}}
 = (z^{n+b}\cdot f_{0}(z))^{\frac{l}{k}};\hspace{.1in} f_{0}(0)\neq 0 \]
implies that $w_{2}$ is a multi--valued function of $z$. Subdivide the
disc $D$ into equal sectors $Q_{m} \ , \ 0\leq m\leq n-1$, and hence define
on each $Q_{m}$ a holomorphic function $F_{m}(z)$, such that $F_{0}$ is 
the principal branch of $(z^{n+b}\cdot f_{0}(z))^{\frac{l}{k}}$, and
$F_{m+1}(z) := e^{2\pi{\bold i}\frac{l}{k}}\cdot F_{m}(z)$. 

As in the previous example, letting $\vartheta = 2\pi t$, we have
\[\omega = \frac{1}{2}(f_{1}+{\bold i}f_{2})(\nu,\bar{\nu})d\bar{\nu} = 
\frac{1-p\cdot\varphi(|\nu|)^{2}}{4\pi l\cdot|\nu|^{2}}\cdot{\bold i}
\nu d\bar{\nu} = \bar{\partial}g(\nu,\bar{\nu}) \ .\]
Let 
\[\gamma(s) := \frac{1-p\cdot\varphi(\sqrt{s})^{2}}{l\cdot s} \ ,\]
noting $\varphi(0) = \frac{1}{\sqrt{p}}$ implies that the improper integral
\[\hat{\gamma}(s) := \int_{0}^{s}\gamma(\tau)d\tau \]
is convergent. It follows that we can set $g(\nu,\bar{\nu}) = 
\frac{{\bold i}}{4\pi}\hat{\gamma}(|\nu|^{2})$. Hence
\[\Phi^{*}\omega = \bar{\partial}(g\circ\Phi) = \frac{{\bold i}}{4\pi}
\bar{\partial}\hat{\gamma}(|\Phi|^{2})\]
is smoothly defined on $D$. Moreover
\[\hat{G} = 2\pi{\bold i}g\circ\Phi = -\frac{1}{2}\hat{\gamma}(|\Phi|
^{2}) \ .\]
Setting $t = \frac{1}{2\pi}arg(z^{n})$, we then have
\[t+{\bold i}\hat{a} = \frac{1}{2\pi{\bold i}}(\log(z^{n}) - \hat{G})
 = \frac{1}{2\pi}arg(z^{n}) +\frac{{\bold i}}{2\pi}(\hat{G} - n\ln|z|) \ ,\]
with which we combine the statement of Theorem 1 to conclude 
\[a(z) = -n\ln(|z|e^{\frac{1}{2n}\hat{\gamma}(|\Phi|^{2})})\]  
in order to define a pseudoholomorphic map $(\psi, a):D\setminus\{0\}
@>>>{\Bbb E}_{p,q}\times{\Bbb R}$ of charge $n$ at the origin, 
asymptotic to the orbit corresponding to $\{w_{2} = 0\}$.
\vspace{.2in}
 
\section{references}

[1] Bourgeois, F. {\em PhD Thesis} (preprint) Stanford 2004

[2] Brieskorn, E. and Kn\"{o}rrer, H. {\em Plane Algebraic Curves}, 
Birkh\"{a}user 1986

[3] Eliashberg, Y. {\em Invariants in Contact Topology}, Proceedings ICM,
Berlin 1998, Volume II, Documenta Mathematica (1998) 327--338

[4] Gluck. H. and Warner, F., {\em Great Circle Fibrations of the 
Three--Sphere}, Duke Math. J. 50 No.1 (1983) 107--132

[5] Hale, J. {\em Ordinary Differential Equations}, Pure and Appl. Math. 21,
Wiley 1969

[6] Hofer, H., {\em Pseudoholomorphic Curves in Symplectisations with 
applications to the Weinstein Conjecture in dimension three}, Invent. Math.
114 (1993) 515--563

[7] Hofer, H. and Kriener, M. {\em Holomorphic Curves in Contact Dynamics},
Proc. Symp. Pure Math. 65 (1999) 77--131

[8] Hofer, H., Wysocki, K. and Zehnder, E. {\em Properties of 
Pseudoholomorphic Curves in Symplectisations I: Asymptotics}, Ann. Inst. 
Henri Poincar\'{e} 13 (1996) 337--379

[9] Hofer, H., Wysocki, K. and Zehnder, E. {\em Properties of 
Pseudoholomorphic Curves in Symplectisations II: Embedding controls and 
algebraic invariants}, Geom. Funct. Anal. 5 (1995) 270--328

[10] Hofer, H., Wysocki, K. and Zehnder, E. {\em Properties of 
Pseudoholomorphic Curves in Symplectisations IV: Asymptotics with 
degeneracies}, in ``Contact and Symplectic Geometry'', C.B Thomas ed.,
Cambridge (1996) 78--117

[11] Martinet, J. {\em Formes de Contact sur les vari\'{e}t\'{e}s de
dimension 3}, Springer Lecture Notes 207 (1971) 142--163

[12] McDuff, D. {\em Singularities of $J$--holomorphic Curves in almost
complex 4--manifolds}, J. Geom. Anal. 2 (1992) 249--266

[13] McKay, B. {\em Dual Curves and Pseudoholomorphic Curves}, Selecta
Math. 9 (2003) 251--311

[14] Micallef, M. and White, B. {\em The structure of branch points in
Minimal Surfaces and in Pseudoholomorphic Curves}, Ann. Math. 139 (1994)
35--85

[15] Sikorav, J.--C. {\em Singularities of $J$--holomorphic Curves}, Math. Z.
226 (1997) 359--373     

\vspace{.1in}

School of Mathematics, Statistics and Computer Science \newline
University of New England \newline
Armidale, NSW 2351 \newline
Australia \newline
adamh@@turing.une.edu.au

\vspace{.1in}

School of Mathematics and Statistics \newline
Melbourne University \newline
Parkville, VIC 3010 \newline
Australia \newline
wysocki@@ms.unimelb.edu.au

\end{document}